\documentclass[a4paper]{article}
\usepackage [a4paper,left=3 cm,bottom=2.5cm,right=3cm,top=2.5cm]{geometry}
\usepackage[english]{babel}
\usepackage{amssymb}
\usepackage{amsmath,amsthm}
\usepackage{subcaption}
\usepackage{tabularx,array}
\usepackage{ulem}
\usepackage{graphicx, url}
\usepackage{enumitem} 
\usepackage{xcolor}

\newcommand{\E}{\mathbb{E}}
\newcommand{\K}{\mathbb{K}}

\newcommand{\R}{\mathbb{R}}

\makeatletter
\newcommand*{\barfix}[2][.175ex]{%
  \mathpalette{\@barfix{#1}}{#2}%
}
\newcommand*{\@barfix}[3]{%
  \vbox{%
    \kern#1\relax
    \hbox{$#2#3\m@th$}%
  }%
}
\makeatother

\newcommand{\modar}{\color{black}}
\newcommand{\modarn}{\color{black}}
\newcommand{\abs}[1]{\left|#1\right|}
\newcommand{\norm}[1]{\left\lVert#1\right\rVert}

\newtheorem{theorem}{Theorem}
\newtheorem{lemma}{Lemma}

\newtheorem{proposition}{Proposition}
\newtheorem{definition}{Definition}

\newcommand{\keywords}[1]{\par\addvspace\baselineskip
\noindent\enspace\ignorespaces#1}

\newcommand{\rev}{\color{black}}
\newcommand{\rew}{\color{black}}
\begin{document}

\title{Estimation of the invariant density for discretely observed diffusion processes: impact of the sampling and of the asynchronicity}
\author{Chiara Amorino\thanks{ Universit\'e du Luxembourg, L-4364 Esch-Sur-Alzette, Luxembourg.The author gratefully acknowledges financial support of ERC Consolidator Grant 815703 “STAMFORD: Statistical Methods for High Dimensional Diffusions”.} \qquad Arnaud Gloter \thanks{Laboratoire de Math\'ematiques et Mod\'elisation d'Evry, CNRS, Univ Evry, Universit\'e Paris-Saclay, 91037, Evry, France.}} 

\maketitle

\begin{abstract}
We aim at estimating in a non-parametric way the density $\pi$ of the stationary distribution of a $d$-dimensional stochastic differential equation $(X_t)_{t \in [0, T]}$, for $d \ge 2$, from the discrete observations of a finite sample $X_{t_0}$, ... , $X_{t_n}$ with $0= t_0 < t_1 < ... < t_n =: T_n$. We propose a kernel density {\rev estimator } and we study its convergence rates for the pointwise estimation of the invariant density under
anisotropic {\rev H\"older} smoothness constraints. First of all, we find some conditions on the discretization step that ensures it is possible to recover the same rates as {\rev if} the continuous trajectory of the process was available. As proven in the recent work \cite{Companion}, such rates are optimal and new in the context of density estimator. Then we deal with the case where such a condition on the discretization step is not satisfied, which we refer to as intermediate regime. In this new regime we identify the convergence rate for the estimation of the invariant density over anisotropic {\rev H\"older} classes, which is the same convergence rate as for the estimation of a probability density belonging to an anisotropic {\rev H\"older} class, associated to $n$ iid random variables $X_1, ..., X_n$. After that we focus on the asynchronous case, in which each component can be observed {\rev at different time points}. Even if the {\rev asynchronicity} of the observations 
{\modarn complexifies the computation of the variance of the estimator, we are able to find conditions ensuring that 
this variance is comparable to the one of the continuous case. We also exhibit that the non synchronicity of the data introduces additional bias terms in the study of the estimator.}
%We also exhibit that the impact of non synchronicity may lead to
%additional non negligible bias terms.}
%implies some difficulties, we are able to overcome them by considering some combinatorics {\rev (see Section \ref{S: proof asynch} for details)} and by proving some sharper bounds on the variance which allow us to lighten the condition needed on the discretization step.

\keywords{Non-parametric estimation, stationary measure, discrete observation, convergence rate, ergodic diffusion, anisotropic density estimation, asynchronous framework.}

\end{abstract}

\section{Introduction}
In this paper we aim at estimating the invariant density belonging to an anisotropic {\rev H\"older} class starting from the discrete observation of the d-dimensional process
\begin{equation}
X_t= X_0 + \int_0^t b( X_s)ds + \int_0^t a(X_s)dW_s, \quad t \in [0,T], 
\label{eq: model intro}
\end{equation}
for $d\ge2$; with $b : \mathbb{R}^d \rightarrow \mathbb{R}^d$, $a : \mathbb{R}^d \rightarrow \mathbb{R}^d \times \mathbb{R}^d$ and $W = (W_t, t \ge 0)$ a d-dimensional Brownian motion. {\rev In this paper we start considering the synchronous case. It means that the observations over the different directions are all recorded at the same instants $0= t_0 \le t_1 \le ... \le t_n=: T_n$. After that we consider the asynchronous case, where the components of the process $X$ are observed in different point of time. 
We assume in this second case that we dispose of the discrete observations $X_{t_1^l}^l, ..., X_{t_n^l}^l$ for any $l \in \{ 1, ..., d \}$, with $0\le t_1^l < ... < t_n^l \le T_n$ for any $l \in \left \{1, ... , d \right \}$. We also define the discretization step $\Delta_n := \sup_{i = 0, ... , n-1} (t_{i + 1} - t_i)$ and $\Delta_n := \sup_{l =1, ... , d} \quad \sup_{i = 0, ... , n} (t_{i + 1}^l - t_i^l)$ with $t_{n+1}^l=T_n$ in the synchronous and asynchronous framework, respectively. The regime considered is $T_n  \rightarrow \infty$ and $\Delta_n \rightarrow 0$ for $n \rightarrow \infty$. }

The model presented in \eqref{eq: model intro} is interesting from a theoretical point of view and because of its applications in many fields. For example, it has application in biology \cite{194 Iacus} and epidemiology \cite{17 Iacus} as well as in physics \cite{176 Iacus} and mechanics \cite{147 Iacus}. Some other classical examples are neurology \cite{109 Iacus}, mathematical finance \cite{115 Iacus} and economics \cite{26 Iacus}.

 Because of the importance of the model, statistical inference for stochastic differential equations has been widely investigated in a lot of different context: disposing of continuous or discrete observations; on a fixed time interval or on long time intervals; in the parametric or in the non parametric frameworks.

Hence, there have been a big amount of papers on the topic. Among them, we quote Comte et al \cite{Com07}, Dalalyan and Reiss \cite{RD}, Genon-Catalot \cite{Gen90}, Gobet et al \cite{Gob04}, Hoffmann \cite{Hof99}, Kessler \cite{Kessler}, Lar\'edo \cite{Lar90}, Marie and Rosier \cite{Nic} and Yoshida \cite{Yos92}.  
%{\rev \sout{The detailed study on stochastic differential equations leads the way for statistical inference for more complicated models such as stochastic partial differential equations (\cite{Alt20}, \cite{Cia17}),  diffusions with mixed effects (\cite{Pic10}, \cite{Del18}), SDEs driven by L\'evy processes (\cite{Mas19}), diffusions with jumps (\cite{SJS}, \cite{Unbiased}, \cite{Mas07}, \cite{Sch19}), Hawkes processes (\cite{Bac13}, \cite{Dio19}, \cite{Hawkes}) and diffusions with discontinuous coefficients (\cite{Sara}, \cite{Lej_Pig18}, \cite{Lej_Pig20}).}}

In this paper we focus on the non parametric estimation of the invariant density starting from the discrete observations of the stochastic differential equation in \eqref{eq: model intro}. In particular, we will consider at the beginning the case where the data is synchronously available and we will study, then, the case where the observations are given asynchronously. The asynchronous case is extremely important for the applications. {\rev Asset prices, indeed, are generally measured when markets close, even if the closing times may be different across markets. In some cases, such as the United States and Japan, the markets do not have any common open hours, while in some other countries there is a partial overlap in the trading hours. As a consequence, the value of portfolios, value at risk measures and hedge strategies appear distort. Even if the prices are quoted at only slightly different times, the observations are still biased. This has been proven to be more evident in analysis of individual markets where closing prices may be stale (see e.g \cite{SchWil77}, \cite{LoMac90}). \\
In today's global market these problems are relevant. Asynchronous data, indeed, makes more complicated the tasks of financial management as the value of the portfolio is never known at a particular time and so the measures such as value at risks and hedge strategies may be misleading.}

 The treatment of non-synchronous trading effects dates back to Fisher \cite{Fis_66}. For several years researchers focused mainly on the effects that stale quotes have on daily closing prices. Campbell et al. (Chapter 3 of \cite{Cam_97}) provides a survey of this literature. There are several models applied on asynchronous financial market closing prices, see for example Chapters 3 and 4 of \cite{As18} and \cite{As2} for respectively autoregressive and log-Gaussian statistical approaches, whilst \cite{As8} proposes a recurrent neural network {\rev approach}. Some other examples in which the authors deal with asynchronicity are \cite{AsB} and \cite{Peloso}.  {\rev  In \cite{BurEng} the authors propose to synchronize prices by computing estimates of the values of assets even when markets are closed, starting from the information given from markets which are open. Correlation and volatility of asynchronous data are also studied in \cite{HayYos05}, where the authors consider the problem of estimating the covariance of two diffusion processes when they are observed at discrete times in a non-synchronous manner. They propose a new estimator which allows them to correct the Epps effect.}

In this context, we aim at proposing a kernel density estimator based on the discrete observations of \eqref{eq: model intro} and we aim at finding the convergence rates of estimation for the stationary measure $\pi$ associated to it, assuming that it belongs to an anisotropic {\rev H\"older} class. As the smoothness properties of elements of a function space may depend on the chosen direction of $\R^d$, the notion of anisotropy plays an important role. We will present some conditions that the discretization step needs to satisfy in order to recover the same (optimal) convergence rate achievable when a continuous trajectory of the process $X$ {\rev is} available and we will discuss the convergence rates we find in the intermediate regime, which is when the discretization step is not small enough and so the discretization error is not negligible. 

With regard to the literature already existing about the estimation of the invariant measure, it is important to say that it is a problem already widely faced in many different frameworks by many authors (see for example \cite{Strauch new}, \cite{Banon}, \cite{Bos98}, \cite{Del80}, \cite{Kut_1998}, \cite{Kut_2004}, \cite{Kut}, \cite{Ngu79} and \cite{Zan01}).
The reason why such a problem results very attractive is the huge amount of physical applications and numerical methods connected to it, such as the Markov Chain Monte Carlo method. For example, {\rev stability theory is used to study properties of the invariant distribution and mixing properties of the diffusion in \cite{Has80} and \cite{Banon}.}
 In \cite{LamPag02} and \cite{Pan08}, instead, it is possible to find an approximation algorithm for the computation of the invariant density. The non-parametric estimation of the invariant density can also be used in order to estimate the drift coefficient in a non-parametric way (see \cite{Nic} and \cite{Sch}). \\ 
Kernel estimators are widely employed as powerful tools: in \cite{Banon} and \cite{Bosq9} some kernel estimators are used to estimate the density of a continuous time process. They are used also in more complicated {\rev models}, such as in jump-diffusion framework (see \cite{Lower bound}, \cite{Optimal}, \cite{Chapitre 4} and \cite{Man10}).

Some references in a context closer to ours are \cite{RD}, \cite{Strauch} and \cite{Companion}. In all three papers kernel density estimators have been used for the study of the convergence rate for the estimation of the invariant density associated to a reversible diffusion process with unit diffusion part (in the first two works) or to the same stochastic differential equation as in \eqref{eq: model intro} (in the third one). They are all based on the continuous observation of the process considered. \\
In particular in \cite{Companion} the invariant density $\pi$ has been estimated by means of the kernel estimator $\hat{\pi}_{h,T}$ assuming to have the continuous record of the process $X$ solution to \eqref{eq: model intro} up to time $T$ and the following upper bound for the mean squared error has been {\rev shown}, for $d \ge 3$:
\begin{equation*}
\sup_{(a,b)\in \Sigma} \mathbb{E} [|\hat{\pi}_{h,T}(x) - \pi (x)|^2] \underset{\sim}{<}
\begin{cases}
(\frac{\log T}{T})^{\frac{2\bar{\beta}_3}{2\bar{\beta}_3 + d - 2}} \qquad \mbox{for } \beta_2 < \beta_3, \\
(\frac{1}{T})^{ \frac{2\bar{\beta}_3}{2\bar{\beta}_3+ d - 2}} \qquad \mbox{for } \beta_2 = \beta_3,
\end{cases}
\end{equation*}
where $\Sigma$ is a class of coefficients for which the stationarity density has some prescribed regularity, $\beta_1 \le \beta_2 \le ... \le \beta_d$ and $\bar{\beta}_3$ is the harmonic mean over the smoothness after having removed the two smallest. In particular, it is $\frac{1}{ \barfix{\bar{\beta}_3}} := \frac{1}{d -2} \sum_{l \ge 3} \frac{1}{\beta_l}$ and so it clearly follows that $\bar{\beta}_3 \ge \bar{\beta}$, where $\frac{1}{\bar{\beta}} := \frac{1}{d} \sum_{j= 1}^d \frac{1}{\beta_j}$. It has also been proven that the convergence rates here above are optimal. 

In this paper, we propose to estimate the invariant density $\pi$ starting from the discrete observations of the process $X$ by means of the kernel estimator $\hat{\pi}_{h,n}$, which is the discretized version of $\hat{\pi}_{h,T}$ (see Section \ref{s: est}).
Then, we prove an upper bound on the variance which is composed of two terms. The first is the same as when the continuous trajectory of the process {\rev is} available, while the second is the discretization error. {\rev If the second is negligible we get a condition the discretization step has to satisfy to recover the continuous convergence rate. On the other side, if the first is negligible compared to the second, we obtain the new convergence rate in the intermediate regime.} In particular, for $d\ge 3$, we show the following: 
\begin{equation*}
\sup_{(a,b)\in \Sigma} \mathbb{E} [|\hat{\pi}_{h,n}(x) - \pi (x)|^2] \underset{\sim}{<}
\begin{cases}
(\frac{\log T_n}{T_n})^{\frac{2\bar{\beta}_3}{2\bar{\beta}_3 + d - 2}} \qquad \mbox{if } \beta_2 < \beta_3 \mbox{ and } \Delta_n \le h_1^* h_2^* \sum_{j = 1}^d |\log h_j^*|, \\
(\frac{1}{T_n})^{ \frac{2\bar{\beta}_3}{2\bar{\beta}_3+ d - 2}} \qquad \mbox{for } \beta_2 = \beta_3 \mbox{ and } \Delta_n \le h_1^* h_2^*,
\end{cases}
\end{equation*}
where $\Delta_n$ is the discretization step and $h^*(T_n) = (h_1^*(T_n), ... , h_d^*(T_n))$ is the rate optimal choice for the bandwidth $h$ (see Theorem \ref{th: discrete d ge 3} and the discussion below for details about the dependence of $h^*$ in $T_n$). \\
On the other side, when the condition above on the discretization step are not respected, we obtain
\begin{equation*}
\sup_{(a,b)\in \Sigma} \mathbb{E} [|\hat{\pi}_{h,n}(x) - \pi (x)|^2] \underset{\sim}{<} n^{- \frac{2\bar{\beta}}{2\bar{\beta}+ d}}.
\end{equation*}
We remark that the convergence rate in the intermediate regime is the same as for the estimation of a probability
density belonging to an anisotropic {\rev H\"older} class, associated to $n$ iid random variables $X_1, ... , X_n$. An analogous result is {\rev shown} also in the bi-dimensional case. \\
\\
After that, we focus on the asynchronous frameworks, in which each component can be observed in a different moment. Such a context implies many difficulties, as in this way even the choice of the estimator to propose appears challenging. The idea is to introduce $d$ functions $\varphi_{n,l} : [0, T] \rightarrow \mathbb{R}$ such that
$$\varphi_{n,l} (t) = \sup \left \{ t_i^l \, | \, t_i^l \le t \right \},$$
 where $(t_i^l)_i$ are the instants of time in which the component $X^l$ is observed, for $l \in \left \{ 1, ... , d \right \}$. In this way it is possible to write the sums as integrals and to propose an estimator which is the natural adaptation of $\hat{\pi}_{h,T}$, the one considered in the continuous framework in \cite{Companion}. Moreover, the non-synchronicity involves some other challenges in the computation of the upper bound on the variance, {\rev most importantly} {\modarn the fact that the observation times of different components are intertwined complexifies the use of the transition density.} We are able to overcome such issues by considering some combinatorics {\rev (see the proof of Proposition \ref{prop: bound var as} below for details)} and by proving some technical sharp bounds which allows us to show that, in the asynchronous context, the condition $\Delta_n \le h_1^* h_2^* $ is enough to recover the variance obtained in the continuous case, for $\beta_2 < \beta_3$.

{ \rev We analyse then the bias term. Here, a condition on the asynchronicity naturally appears. In particular, having defined 
$$\Delta_n':= \sup_{t \in [0, T_n]} \sup_{i, j=1, ... , d} |\varphi_{n,i}(t) - \varphi_{n,j}(t)|,$$ we obtain the continuous convergence rate when additionally the condition $\Delta_n' \le (\frac{\log T_n}{T_n})^{\frac{2 \bar{\beta}_3}{2 \Bar{\beta}_3 + d - 2}}$ holds. 
%\sout{From both our findings here above and the results in \cite{Companion} it appears clearly that the first two components do not have the same weight as all the others. Hence, we decide to deal with a different framework, in which the first two components are observed continuously and all the others are observed in a discrete way. In this context we are able to remove, under a non restrictive condition on $\Delta$, the discretization error in the bound on the variance and to recover the same convergence rates as when the continuous trajectory of the process was available. We also show that the bias term surprisingly provides a strong condition on the discretization step. We remark that, when the continuous trajectory of the process is available, the variance is the most difficult term to upper bound, while the bound on the bias term is classical. The constraint on the discretization step derives from the asynchronicity of the problem. However, it is enough to ask $\Delta_n \le (\frac{\log T_n}{T_n})^{\frac{2\bar{\beta}_3}{2\bar{\beta}_3 + d - 2}}$ to recover the same convergence rate as when the continuous trajectory of the process is available. }
}

The outline of the paper is the following. In Section \ref{s: model} we introduce the model and we list the assumptions we will need in the sequel, while in Section \ref{S: Continuous observations} we recall the results in the case where the continuous trajectory of the process $X$ is available. In Section
\ref{S: Discrete observations} we consider the synchronous framework. We propose the kernel estimator and we state the upper bounds for the variance which will result in conditions on the discretization step to obtain the continuous regime and in the convergence rates in the intermediate regime. Section \ref{s: asynch} is devoted to the statement of our results in the asynchronous framework. In Section \ref{S: proofs synch} we prove the results stated in Section \ref{S: Discrete observations} while Section \ref{S: proof asynch} is devoted to the proof of results under asynchronicity.

\section{Model Assumptions}{\label{s: model}}
We aim at proposing a non-parametric estimator for the invariant density associated to a {\rev $d$}-dimensional diffusion process $X$ {\rev discretely observed}. In {\rev Section \ref{S: Continuous observations}} we will recall what happens when a continuous record of the process $X^T = \left \{  X_t, 0 \le t \le T \right \}$ up to time $T$ is available. After that, {\rev in Section \ref{S: Discrete observations}}, we will be working in a high frequency setting and we will wonder which conditions on the discretization step will ensure the achievement of the same results as in the continuous case.
The diffusion is a strong solution of the following stochastic differential equation:
\begin{equation}
X_t= X_0 + \int_0^t b( X_s)ds + \int_0^t a(X_s)dW_s, \quad t \in [0,T], 
\label{eq: model}
\end{equation}
where $b : \mathbb{R}^d \rightarrow \mathbb{R}^d$, $a : \mathbb{R}^d \rightarrow \mathbb{R}^d \times \mathbb{R}^d$ and $W = (W_t, t \ge 0)$ is a d-dimensional Brownian motion. The initial condition $X_0$ and $W$ are independent. We denote $\tilde{a}:= a \cdot a^T$. \\

We denote with $|.|$ and $<., . >$ respectively the Euclidian norm and the scalar product in $\mathbb{R}^d$, and for a matrix in $\mathbb{R}^d \otimes \mathbb{R}^d$ we denote its operator norm by  $|.|$. \\ \\
\textbf{A1}: \textit{The functions $b(x)$ and $a(x)$ are bounded globally Lipschitz functions of class $\mathcal{C}^1$, such that for all $x\in\mathbb{R}^d$,
	$$ |a(x)|\le a_0,~ |b(x)|\le b_0,~|\frac{\partial}{\partial x_i} b(x)|\le b_1, ~ |\frac{\partial}{\partial x_i} a(x)|\le a_1 ,\text{ for $i \in \{1,\dots,d\}$},$$
	where $a_0>0,~b_0>0,~a_1>0,~b_1>0$ are some constants.
	% and  $|.|$ is the Euclidian norm.
	Moreover, for some $a_{\text{min}} > 0$,
	$$a_{\text{min}}^2 \mathbb{I}_{d \times d} \le \tilde{a}(x) $$%\le c \mathbb{I}_{d \times d}, $$
	where $\mathbb{I}_{d \times d}$ denotes the $d \times d$ identity matrix.} \\ \\
	{\rev As the inequality here above is between two {\rew matrices}, it is worth explaining it is intended in the scalar product sense.}\\ \\
	\textbf{A2 (Drift condition) }: \textit{ \\
 There exist $\tilde{C}_{\text{b}} > 0$ and $\tilde{\rho}_{\text{b}} > 0$ such that $<x, b(x)>\, \le -\tilde{C}_{\text{b}} |x|$, $\forall x : |x| \ge \tilde{\rho}_{\text{b}}$.
 } \\
\\
We suppose that the invariant probability measure $\mu$ of $X$ is absolutely continuous with respect to the Lebesgue measure and from now on we will denote its density as $\pi$: {\rev $ \mu(dx) = \pi(x) dx$.}
Under the assumptions A1 - A2 the process $X$ admits a unique invariant distribution $\mu$ and the ergodic theorem holds. {\rev In particular, it implies the exponential ergodicity of the process $X$. For the exponential mixing property of general multidimensional diffusions, the reader may consult Theorem 3 of Kusuoka and Yoshida \cite{Kus_Yos} ($\alpha$ - mixing); Meyn and Tweedie \cite{May_Twe}, Stramer and Tweedie \cite{Str_Twe} and Veretennikov \cite{Ver} ($\beta$ - mixing). \\ %We recall below the notion of exponential ergodicity. \\
{\rew Under our assumptions the process $X$ is exponentially $\beta$ mixing and exponentially ergodic. In particular, the following inequality holds true:
$$\left \| P_{t} f \right \|_{L^1 (\mu)} \le c e^{- \rho t} \left \| f \right \|_{\infty}, $$} %a process $X$ is if it admits a unique invariant distribution $\pi$ and additionally if there exist positive constants c and $\rho$ for which, for each $f$ centered under $\mu$, 
where $P_{t} f(x) := \mathbb{E} [f(X_t) | X_0 = x ]$ is the transition semigroup of the process $X$. \\
The transition density is denoted by $p_{t}$ and it is such that $P_{t} f(x) = \int_{\mathbb{R}^d} f(y) p_{t} (x,y) dy$.}  \\
\\
We want to estimate the invariant density $\pi$ belonging to the anisotropic H{\"o}lder class $\mathcal{H}_d (\beta, \mathcal{L})$ defined below.
\begin{definition}
Let $\beta = (\beta_1, ... , \beta_d)$, $\beta_i > 0$, $\mathcal{L} =(\mathcal{L}_1, ... , \mathcal{L}_d)$, $\mathcal{L}_i > 0$. A function $g : \mathbb{R}^d \rightarrow \mathbb{R}$ is said to belong to the anisotropic H{\"o}lder class $\mathcal{H}_d (\beta, \mathcal{L})$ of functions if, for all $i \in \left \{ 1, ... , d \right \}$,
$$\left \| D_i^k g \right \|_\infty \le \mathcal{L}_i \qquad \forall k = 0,1, ... , \lfloor \beta_i \rfloor, $$
$$\left \| D_i^{\lfloor \beta_i \rfloor} g(. + t e_i) - D_i^{\lfloor \beta_i \rfloor} g(.) \right \|_\infty \le \mathcal{L}_i |t|^{\beta_i - \lfloor \beta_i \rfloor} \qquad \forall t \in \mathbb{R},$$
for $D_i^k g$ denoting the $k$-th order partial derivative of $g$ with respect to the $i$-th component, $\lfloor \beta_i \rfloor$ denoting the largest integer strictly smaller than $\beta_i$ and $e_1, ... , e_d$ denoting the canonical basis in $\mathbb{R}^d$.
\end{definition}

	This leads us to consider a class of coefficients $(a,b)$ for which the stationary density $\pi=\pi_{(a,b)}$ has some prescribed H\"older regularity.
 	\begin{definition}\label{D:def_Sigma}
		Let $\beta = (\beta_1, ... , \beta_d)$, $0 < \beta_1 \le ... \le \beta_d$ and $\mathcal{L} =(\mathcal{L}_1, ... , \mathcal{L}_d)$, $\mathcal{L}_i > 0$, $0<a_\text{min}\le a_0$ and $a_1>0, ~ b_0>0, ~b_1>0, ~\tilde{C}>0, ~\tilde{\rho}>0$.
		
		We define $\Sigma (\beta, \mathcal{L},a_{\text{min}},a_0,a_1,b_0,b_1,\tilde{C},\tilde{\rho})$ the set of couple of functions $(a,b)$ where  $a: \mathbb{R}^d \rightarrow \mathbb{R}^d\otimes\mathbb{R}^d$ and $b: \mathbb{R}^d \rightarrow \mathbb{R}^d$ are such that 
		\begin{itemize}
			\item 		$a$ and $b$ satisfy A1 with the constants $(a_\text{min},a_0,a_1,b_0,b_1)$,
			\item 		$b$ satisfies A2 with the constants $(\tilde{C},\tilde{\rho})$,
			\item  the density $\pi_{(a,b)}$ of the invariant measure associated to the stochastic differential equation \eqref{eq: model} belongs to $\mathcal{H}_d (\beta, \mathcal{L})$.
		\end{itemize}
		\label{def: insieme sigma v2}
	\end{definition}

{\rev To have an idea about the link between the coefficients and the invariant density, one can think about a reversible diffusion. Indeed, if the diffusion coefficient is the identical matrix and the drift is such that $b(x) = - \nabla V(x)$, where $V$ is a function we refer to as a potential, then assuming $V \in \mathcal{H}_d (\beta, \mathcal{L})$ implies ${ \pi \modarn = C e^{-V} } \in \mathcal{H}_d (\beta + 1, \mathcal{L})= \mathcal{H}_d (\beta_1 + 1, ... , \beta_d + 1, \mathcal{L})$ .}\\
\\
We aim at estimating the invariant density $\pi$ starting from discrete observations of the process $X$. In particular, we want to find some conditions that the discretization step has to satisfy in order to recover the same convergence rates we had when a continuous record of the process was available. Moreover, one may wonder which are the convergence rates in intermediate regime, i.e. when the discretization step goes to zero but the associated error is not negligible. In order to answer these questions we recall what happens when the whole trajectory of the process $X$ is available, as detailed discussed in \cite{Companion}. This is the purpose of {\rev the} next section.

\section{Continuous observations}\label{S: Continuous observations}

 It is natural to estimate the invariant density $\pi \in \mathcal{H}_d (\beta, \mathcal{L})$ by means of a kernel estimator.
We therefore introduce some kernel function $K: \mathbb{R} \rightarrow \mathbb{R}$ satisfying 
\begin{equation}
\int_\mathbb{R} K(x) dx = 1, \quad \left \| K \right \|_\infty < \infty, \quad \mbox{supp}(K) \subset [-1, 1], \quad \int_\mathbb{R} K(x) x^l dx = 0,
\label{eq: properties K}
\end{equation}
for all $l \in \left \{ 0, ... , M \right \}$ with $M \ge \max_i \beta_i$. \\
For $j \in \left \{ 1, ... , d \right \}$, we denote by $X_t^j$ the $j$-th component of $X_t$. A natural estimator of $\pi$ at $x= (x_1, ... , x_d)^T \in \mathbb{R}^d$ in the anisotropic context is given by 
\begin{equation}
\hat{\pi}_{h,T}(x) = \frac{1}{T \prod_{l = 1}^d h_l} \int_0^T \prod_{m = 1}^d K(\frac{x_m - X_u^m}{h_m}) du.
\label{eq: def estimator}
\end{equation}
The multi-index $h=(h_1, ... , h_d) $ is small. In particular, we assume $h_i < 1$ for any $i \in \{1, ... , d \} $. \\
\\
For $d \ge 3$, from Theorem 1 of \cite{Companion} we have the following convergence rate for the kernel estimator proposed in \eqref{eq: def estimator} and for the optimal bandwidth given below, in \eqref{eq: optimal bandwidth}: 
\begin{equation*}
 \sup_{(a, b) \in \Sigma} \mathbb{E}[|\hat{\pi}_{h,T}(x) - \pi (x)|^2] \underset{\sim}{<} 
 \begin{cases}
 (\frac{\log T}{ T})^ {\frac{2\bar{\beta}_3}{2\bar{\beta}_3 + d - 2}} \qquad \mbox{if } \beta_2 < \beta_3 \\
 (\frac{1}{ T})^ {\frac{2\bar{\beta}_3}{2\bar{\beta}_3 + d - 2}} \qquad \mbox{if } \beta_2 = \beta_3, 
 \end{cases}
\end{equation*}
where $\Sigma$ is the set defined in Definition \ref{def: insieme sigma v2} and $\bar{\beta}_3$ is such that
$$\frac{1}{\bar{\beta}_3} := \frac{1}{d-2} \sum_{j = 3}^d \frac{1}{\beta_j}.$$
Moreover, from Theorems 3 and 4 of \cite{Companion} we know they are optimal {\rev in the minimax sense}. \\
\\
Regarding the bi-dimensional case we know, from Theorem 2 in \cite{Companion} that the following holds true
$$ \sup_{(a, b) \in \Sigma} \mathbb{E}[|\hat{\pi}_{h,T}(x) - \pi (x)|^2] \underset{\sim}{<} \frac{\log T}{ T}$$
and the convergence rate here above is optimal {\rev in the minimax sense} (see Theorem 5 of \cite{Companion}).\\
\\
  We now suppose that the continuous record of the process, up to time $T$, is no longer available. {\rev We dispose instead} of the discretization of the process at the instants $0= t_0 \le t_1 \le ... \le t_n=T$, given by $X_{t_0}, ..., X_{t_n}$. 
The first goal of Section \ref{S: Discrete observations} is to find some conditions the discretization step has to satisfy in order to recover the same convergence rates as in this section.

\section{Discrete observations, synchronous framework}{\label{S: Discrete observations}}
In this section we suppose that we observe a finite sample $X_{t_0}, ..., X_{t_n}$, with $0= t_0 \le t_1 \le ... \le t_n=: T_n$. The process $X$ is {\rev the} solution of the stochastic differential equation \eqref{eq: model}. Every observation time point depends also on $n$ but, in order to simplify the notation, we suppress this index. We assume the discretization scheme to be uniform which means that, for any $i \in \{ 0, ..., n-1 \}$, it is $t_{i + 1} - t_i = : \Delta_n$. We will be working in a high-frequency setting i.e. the discretization step $\Delta_n \rightarrow 0$ for $n \rightarrow \infty$. We assume moreover that $T_n = n \Delta_n \rightarrow \infty$ for $n \rightarrow \infty$ and that $\Delta_n > n^{- k}$ for some $k \in(0, 1)$.  

\subsection{Construction estimator}{\label{s: est}}
As in Section \ref{S: Continuous observations}, we propose to estimate the invariant density $\pi \in \mathcal{H}_d (\beta, \mathcal{L})$ associated to the process $X$, solution to \eqref{eq: model}. To do that, we propose a kernel estimator which is the discretized version of the one introduced in \eqref{eq: def estimator}. For $x=(x_1, ..., x_d) \in \R^d$, we define 
\begin{align}\label{eq: def discrete estimator} 
\hat{\pi}_{h, n} (x) & := \frac{1}{n \Delta_n} \frac{1}{\prod_{l = 1}^d h_l} \sum_{i = 0}^{n-1}  \prod_{l = 1}^d K(\frac{x_l - X_{t_i}^l}{h_l}) (t_{i + 1} - t_i) \\
&= \frac{1}{n} \sum_{i = 0}^{n-1} \mathbb{K}_h(x - X_{t_i}), \nonumber
\end{align}
with $K$ a kernel function as in \eqref{eq: properties K} 
{\modarn and $\mathbb{K}_h=\prod_{l=1}^{d} K_{h_l}$ where $K_{h_l}(\cdot)=\frac{1}{h_l}K(\frac{\cdot}{h_l})$.}   \\
\\
In this context we have two objectives:
\begin{enumerate}
    \item Find some conditions on $\Delta_n$ to get the same convergence rates we had when a continuous record of the process was available.
    \item Find the convergence rates in the intermediate regime ($\Delta_n$ tends to zero but the discretization error is not negligible).
\end{enumerate}
To achieve them, we have to study the behaviour of the mean squared error $\mathbb{E} [|\hat{\pi}_{n,h}(x) - \pi (x)|^2] $ when $d \ge 3$ and when $d=2$, and we have deal with two asymptotic regimes. 
The idea mainly consists in some upper bounds for the variance of the estimator.
The difference, with respect to the continuous case, is that now we get an extra term which derives from the discretization. \\
If the new discretization term is negligible compared to the others, then the convergence rates are the same they were in the continuous case; otherwise we will find new convergence rates.

\subsection{Main results, synchronous framework}{\label{S: main}}
The asymptotic behaviour of the estimator proposed in \eqref{eq: def discrete estimator} is based on the bias-variance decomposition. To find the convergence rates it achieves we need a bound on the variance, as stated in the following propositions. We recall that the estimator performs differently depending on the dimension $d$. In this paper we provide the main results for $d \ge 2$. 

\begin{proposition}
Suppose that A1-A2 hold with some constant $0<a_\text{min}\le a_0$ and $a_1>0, ~ b_0>0, b_1>0, \tilde{C}_b, \tilde{\rho}_b$ and that $d \ge 3$. {\rew Let $\beta = (\beta_1, ... , \beta_d)$, $\beta_1 = \beta_2 = ... = \beta_{k_0} < \beta_{k_0 + 1} \le ... \le \beta_d$, for some $k_0 \in \{ 1, ..., d\}$ and $\mathcal{L} =(\mathcal{L}_1, ... , \mathcal{L}_d)$, $\mathcal{L}_i > 0$.}
 If $\hat{\pi}_{h,n}$ is the estimator proposed in \eqref{eq: def discrete estimator}, then there exist $c > 0$ and $T_0 > 0$ such that, for $T_n \ge T_0$, the following holds true.
\begin{itemize}
    \item[$\bullet$] If $k_0 = 1$ and $\beta_2 < \beta_3$ or $k_0 = 2$, then 
  $$Var(\hat{\pi}_{h,n}(x)) \le \frac{c}{T_n} \frac{\sum_{j = 1}^d |\log(h_j)|}{\prod_{l = 3}^d h_l} + \frac{c}{T_n} \frac{\Delta_n}{\prod_{l = 1}^d h_l}.$$   
    \item[$\bullet$] If $k_0 \ge 3$, then 
  $$Var(\hat{\pi}_{h,n}(x)) \le \frac{c}{T_n} \frac{1}{(\prod_{l = 1}^{k_0} h_l)^{1 - \frac{2}{k_0}}(\prod_{l \ge k_0 + 1} h_l)} + \frac{c}{T_n} \frac{\Delta_n}{\prod_{l = 1}^d h_l}.$$  
  \item[$\bullet$] If otherwise $k_0 = 1$ and $\beta_2 = \beta_3$, then 
$$Var(\hat{\pi}_{h,T}(x)) \le \frac{c}{T_n} \frac{1}{\prod_{l \ge 4} h_l \sqrt{h_2 h_3}} + \frac{c}{T_n} \frac{\Delta_n}{\prod_{l = 1}^d h_l}.$$
Moreover, the constant $c$ is uniform over the set of coefficients $(a, b) \in \Sigma$.
\end{itemize}
\label{prop: bound var discrete dge3}
\end{proposition}

{\rev Comparing the result here above with Proposition 2 of \cite{Companion}, it is easy to see that we have an extra term coming from the discretization: $\frac{c}{T_n} \frac{\Delta_n}{\prod_{l = 1}^d h_l}$. \\
One can remark that, as in the high frequency framework the estimators can be understood as numerical approximators of the the full observation estimators, a natural choice could be to derive the results in the discrete case from the continuous one. {\rew However, following the same route as Theorem 4.1 in \cite{Strauch} to control the error coming from the approximation of the estimator, one can check that this approach leads us to a larger discretization error, which yields a worse condition on the discretization step. }
This justifies our choice of proving some upper bounds directly on the variance of the discrete estimator, as in Proposition \ref{prop: bound var discrete dge3}.} Such bounds leads us to the first main result of this section.
\begin{theorem}{[Discretization term negligible] \\}
Suppose that A1-A2 hold and that $d\ge 3$. 
%{\rew Let $\beta = (\beta_1, ... , \beta_d)$, $0 < \beta_1 \le ... \le \beta_d$ and $\mathcal{L} =(\mathcal{L}_1, ... , \mathcal{L}_d)$, $\mathcal{L}_i > 0$.}
{\rew Let $\beta = (\beta_1, ... , \beta_d)$, $0 < \beta_1 \le ... \le \beta_d$, $\mathcal{L} =(\mathcal{L}_1, ... , \mathcal{L}_d)$, $\mathcal{L}_i > 0$, $0<a_\text{min}\le a_0$,
	$a_1>0, ~ b_0>0, ~b_1>0, ~\tilde{C}>0, ~\tilde{\rho}>0$ and
	set $\Sigma=\Sigma (\beta, \mathcal{L},a_{\text{min}},a_0,a_1,b_0,b_1,\tilde{C},\tilde{\rho})$, using the notation of Definition \ref{D:def_Sigma}.} 
Let $h^* = (h_1^*, ... , h_d^*)$ be the rate optimal choice for the bandwidth $h$ as given in \eqref{eq: optimal bandwidth}.
 Then, there exist $c > 0$ and $T_0 > 0$ such that, for $T_n \ge T_0$, the following hold true.
\begin{itemize}
\item[$\bullet$] If $\beta_2 < \beta_3$ and $\Delta_n \underset{\sim}{<} h_1^* h_2^* \sum_{j = 1}^d |\log h_j^*|$, then  
    $$\sup_{(a, b) \in \Sigma}\mathbb{E}[|\hat{\pi}_{h,n}(x) - \pi (x)|^2] {\rew \le c\, }
    (\log T_n / T_n)^{ \frac{2\bar{\beta}_3}{2\bar{\beta}_3 + d - 2}}. $$
    \item[$\bullet$] If otherwise $\beta_2 = \beta_3$ and $\Delta_n \underset{\sim}{<} h_1^* h_2^*$, then 
$$ \sup_{(a, b) \in \Sigma} \mathbb{E}[|\hat{\pi}_{h^*,n}(x) - \pi (x)|^2] {\rew \le c \,}
T_n^{ - \frac{2\bar{\beta}_3}{2\bar{\beta}_3 + d - 2}},$$
{\rev where $\bar{\beta}_3$ is such that $\frac{1}{\bar{\beta}_3} = \frac{1}{d - 2} \sum_{l = 3}^d \frac{1}{\beta_l}$.}
\end{itemize}
\label{th: discrete d ge 3}
\end{theorem}

Comparing the results here above with the ones included in Section 3 of \cite{Companion} we deduce that, when the discretization step satisfies the constraint $\Delta_n \underset{\sim}{<} h_1^* h_2^* \sum_{j = 1}^d |\log h_j^*|$ (or $\Delta_n \underset{\sim}{<} h_1^* h_2^*$ respectively) it is possible to recover the same optimal convergence rates we had when the trajectory of the process was observed continuously.  \\
As seen in the proof of Theorem 1 of \cite{Companion}, for $\beta_2 < \beta_3$ the rate optimal choice for the bandwidth $h$ is given by $h^*_j = (\frac{\log T_n}{T_n})^{a_j}$, while for $\beta=\beta_3$ it is $h^*_j = (\frac{1}{T_n})^{a_j}$ with 
\begin{equation}
a_j = \frac{\bar{\beta}_3}{\beta_j(2 \bar{\beta}_3 + d - 2)} \qquad \mbox{for any } j\in \{1, ..., d\}.
\label{eq: optimal bandwidth}
\end{equation}
We remark that, according to \cite{Companion}, it is also possible {\modarn choose smaller bandwidths $h_1^*$ and $h_2^*$} in the case $\beta_2 < \beta_3$. However, in order to make the condition on the discretization step as weak as possible, it is convenient to choose $h_1^* h_2^*$ as large as possible, which leads us to the choice gathered in \eqref{eq: optimal bandwidth}. \\
Hence, replacing the optimal choice for the bandwidth as in \eqref{eq: optimal bandwidth} one can recover the same upper bound for the mean squared error as in Theorem 1 of \cite{Companion} when the following conditions hold: 
\begin{align}{\label{eq: cond delta T}}
\Delta_n \underset{\sim}{<} (\frac{\log T_n}{T_n})^{\frac{\bar{\beta}_3}{2 \bar{\beta}_3 + d - 2}(\frac{1}{\beta_1} + \frac{1}{\beta_2})} \log T_n \qquad \mbox{for } \beta_2 < \beta_3, \\
\Delta_n \underset{\sim}{<} (\frac{1}{T_n})^{\frac{\bar{\beta}_3}{2 \bar{\beta}_3 + d - 2}(\frac{1}{\beta_1} + \frac{1}{\beta_2})}\qquad \mbox{for } \beta_2 = \beta_3.  {\label{eq: cond delta T 2}}
\end{align}

When such conditions are not respected, instead, we get a different convergence rate. It is achieved by choosing the optimal bandwidth as $h_j(n):= (\frac{1}{n})^{\frac{\bar{\beta}}{\beta_j (2 \bar{\beta} + d)}}$, where $\bar{\beta}$ is the harmonic mean over the $d$ different smoothness: 
$$\frac{1}{ \bar{\beta}} = \sum_{j = 1}^d \frac{1}{\beta_j}. $$
It leads us to the following result {\rew which can be applied when the conditions on the discretization step gathered in Theorem \ref{th: discrete d ge 3} are not respected.}
\begin{theorem}{[Discretization term non-negligible] \\}
%Suppose that A1-A2 hold and 
Suppose that $d \ge 3$. 
%{\rew Let $\beta = (\beta_1, ... , \beta_d)$, $0 < \beta_1 \le ... \le \beta_d$ and $\mathcal{L} =(\mathcal{L}_1, ... , \mathcal{L}_d)$, $\mathcal{L}_i > 0$.} 
{\rew Let $\beta = (\beta_1, ... , \beta_d)$, $0 < \beta_1 \le ... \le \beta_d$, $\mathcal{L} =(\mathcal{L}_1, ... , \mathcal{L}_d)$, $\mathcal{L}_i > 0$, $0<a_\text{min}\le a_0$
	$a_1>0, ~ b_0>0, ~b_1>0, ~\tilde{C}>0, ~\tilde{\rho}>0$ and
	set $\Sigma=\Sigma (\beta, \mathcal{L},a_{\text{min}},a_0,a_1,b_0,b_1,\tilde{C},\tilde{\rho})$, using the notation of Definition \ref{D:def_Sigma}.} 
{\rew Assume that} 
%that the conditions on the discretization step gathered in Theorem \ref{th: discrete d ge 3} are not respected, and  
%so 
one of the following holds
\begin{itemize}
\item[$\bullet$] $\Delta_n > (\frac{\log T_n}{T_n})^{\frac{\bar{\beta}_3}{2 \bar{\beta}_3 + d -2}(\frac{1}{\beta_1} + \frac{1}{\beta_2})} \log T_n$ and $\beta_2 < \beta_3$.
\item[$\bullet$] $\Delta_n > (\frac{1}{T_n})^{\frac{\bar{\beta}_3}{2 \bar{\beta}_3 + d -2}(\frac{1}{\beta_1} + \frac{1}{\beta_2})}$ and $\beta_2 = \beta_3$.
\end{itemize}
 Then, there exist $c > 0$ and $T_0 > 0$ such that, for $T_n \ge T_0$,
$$\sup_{(a,b) \in \Sigma} \mathbb{E}[|\hat{\pi}_{h,n}(x) - \pi (x)|^2] {\rew \le c \,} n^{-\frac{2 \bar{\beta}}{2 \bar{\beta} + d}},$$
where $\bar{\beta}$ is the harmonic mean of the smoothness over the $d$ direction, defined as
$$\frac{1}{\bar{\beta}}:= \frac{1}{d} \sum_{j = 1}^d \frac{1}{\beta_j}.$$
\label{th: discrete d ge 3 non negl}
\end{theorem}

It is interesting to remark that it is also the convergence rate for the estimation of a probability density belonging to an {\rev H\"older} class, associated to $n$ independent and identically distributed random variables $X_1, ..., X_n$. \\
\\
{\rev One may wonder if the rate obtained in the intermediate regime is optimal. The question is addressed for the case $d=1$ in \cite{Malliavin}. Therein it is shown that the convergence rate here above is optimal in a minimax sense. However, the approach used in \cite{Malliavin} relies on local time and on some bounds on the occupation time which are not easily extended to higher dimension. Hence, the question of optimality of the rate obtained when the discretization step is not negligible is still an open question for $d\ge 2$.} \\
\\
{\rev We remark that, as $T_n = n \Delta_n$, the conditions \eqref{eq: cond delta T}, \eqref{eq: cond delta T 2} can be written {\rew as} function of $n$. Using that $\frac{d}{\bar{\beta}} = \frac{d-2}{\bar{\beta}_3} + (\frac{1}{\beta_1} + \frac{1}{\beta_2})$ 
	{\modarn and $\Delta_n >n^{-k}$ with $k \in (0,1)$,}	
	it is possible to check they are equivalent to the following: 
$$\Delta_n \underset{\sim}{<} (\frac{1}{n})^{\frac{\bar{\beta}}{2 \bar{\beta} + d}(\frac{1}{\beta_1} + \frac{1}{\beta_2})} {\modarn \log n } \qquad \mbox{for } \beta_2 < \beta_3,$$
\begin{equation}
\Delta_n \underset{\sim}{<} (\frac{1}{n})^{\frac{\bar{\beta}}{2 \bar{\beta} + d}(\frac{1}{\beta_1} + \frac{1}{\beta_2})} \qquad \mbox{for } \beta_2 = \beta_3.
\label{eq: 8.5}
\end{equation}
We are interested in studying the extreme case, for which the discretization step is equal to the right hand side of the equations here above. Then, the continuous convergence rate always becomes $(\frac{1}{n})^{\frac{2 \bar{\beta}}{2 \bar{\beta} + d}}$. Indeed,
if $\beta_2 = \beta_3$ we have 
$$(\frac{1}{T_n})^{ \frac{2\bar{\beta}_3}{2\bar{\beta}_3 + d - 2}} = (\frac{1}{n \Delta_n})^{ \frac{2\bar{\beta}_3}{2\bar{\beta}_3 + d - 2}} = ((\frac{1}{n})^{1 - \frac{\bar{\beta}}{2 \bar{\beta} + d}(\frac{1}{\beta_1} + \frac{1}{\beta_2})})^{\frac{2\bar{\beta}_3}{2\bar{\beta}_3 + d - 2}} = ((\frac{1}{n})^{\frac{2 + \frac{d-2}{\Bar{\beta}_3}}{2 + \frac{d}{\bar{\beta}}}})^{\frac{2\bar{\beta}_3}{2\bar{\beta}_3 + d - 2}} = (\frac{1}{n})^{\frac{2 \bar{\beta}}{2 \bar{\beta} + d}}.$$
The very same computation holds for $\beta_2 < \beta_3$ as 
$$(\frac{\log T_n}{T_n})^{ \frac{2\bar{\beta}_3}{2\bar{\beta}_3 + d - 2}} = (\log n \Delta_n)^{ \frac{2\bar{\beta}_3}{2\bar{\beta}_3 + d - 2}}((\frac{1}{n})^{1 - \frac{\bar{\beta}}{2 \bar{\beta} + d}(\frac{1}{\beta_1} + \frac{1}{\beta_2})} \frac{1}{\log n \Delta_n} )^{ \frac{2\bar{\beta}_3}{2\bar{\beta}_3 + d - 2}}.$$
It follows that, when the discretization step reaches the threshold, the continuous convergence rate is equal to the convergence rate in the intermediate regime. } \\
\\
 For $d=2$, analogous results hold. In particular, we have the following proposition. 
\begin{proposition}
Suppose that A1-A2 hold and that $d = 2$. {\rew Let $\beta = (\beta_1, ... , \beta_d)$, $0 < \beta_1 \le ... \le \beta_d$ and $\mathcal{L} =(\mathcal{L}_1, ... , \mathcal{L}_d)$, $\mathcal{L}_i > 0$.} Let $\hat{\pi}_{h,n}$ be the estimator proposed in \eqref{eq: def discrete estimator}, then there exist $c > 0$ and $T_0 > 0$ such that, for $T_n \ge T_0$,
$$Var(\hat{\pi}_{h,n}(x)) \le \frac{c}{T_n} \sum_{j = 1}^d |\log(h_j)| + \frac{1}{T_n} \frac{\Delta_n}{h_1 h_2},$$
where the constant $c$ is uniform over the set of coefficients $(a,b) \in \Sigma$.
\label{prop: bound var discrete d=2}
\end{proposition}
As before, it leads us to a condition on $\Delta_n$ that allows us to recover the continuous convergence rate gathered in Theorem 2 of \cite{Companion}, in the continuous case.
\begin{theorem}
%Suppose that A1-A2 hold and that {\rew $d=2$. Let $\beta = (\beta_1, ... , \beta_d)$, $0 < \beta_1 \le ... \le \beta_d$ and $\mathcal{L} =(\mathcal{L}_1, ... , \mathcal{L}_d)$, $\mathcal{L}_i > 0$.} 
{\rew Suppose that $d=2$. Let $\beta = (\beta_1,\beta_2)$, $0 < \beta_1 \le \beta_2$, $\mathcal{L} =(\mathcal{L}_1, \mathcal{L}_2)$, $\mathcal{L}_i > 0$, $0<a_\text{min}\le a_0$
		$a_1>0, ~ b_0>0, ~b_1>0, ~\tilde{C}>0, ~\tilde{\rho}>0$ and
		set $\Sigma=\Sigma (\beta, \mathcal{L},a_{\text{min}},a_0,a_1,b_0,b_1,\tilde{C},\tilde{\rho})$, using the notation of Definition \ref{D:def_Sigma}.} 
	
Let $h^* = (h_1^* , h_2^*)$ be the rate optimal choice \eqref{eq: optimal bandwidth} for the bandwidth $h$. Then, here exist $c > 0$ and $T_0 > 0$ such that, for $T_n \ge T_0$, the following hold true.
\begin{itemize}
    \item[$\bullet$] If $\Delta_n \le h_1^* h_2^* \sum_{j = 1}^2 |\log h_j^*| = (\frac{\log T_n}{T_n})^{\frac{1}{\bar{\beta}}} \log T_n$, then 
    $$\sup_{(a,b) \in \Sigma} \mathbb{E}[|\hat{\pi}_{h,n}(x) - \pi (x)|^2] \le c \frac{\log T_n}{T_n}$$
    \item[$\bullet$] If otherwise $\Delta_n > (\frac{\log T_n}{T_n})^{\frac{1}{\bar{\beta}}} \log T_n $, then 
    $$\sup_{(a,b) \in \Sigma} \mathbb{E}[|\hat{\pi}_{h,n}(x) - \pi (x)|^2] \le c (\frac{1}{n})^{\frac{2 \bar{\beta}}{2 \bar{\beta} + 2}}.$$
\end{itemize}
\label{th: discrete d =2}
\end{theorem}

{\rev As below Theorem \ref{th: discrete d ge 3 non negl}, we can write the conditions on $\Delta$ in function of $n$. In particular, the condition $\Delta_n \le(\frac{\log T_n}{T_n})^{\frac{1}{\bar{\beta}}} \log T_n$ is equivalent to $\Delta_n \le (\frac{1}{n})^{\frac{1}{\bar{\beta} + 1}} \log(n \Delta_n)$. Replacing $\Delta_n = (\frac{1}{n})^{\frac{1}{\bar{\beta} + 1}} \log(n \Delta_n)$ in the continuous convergence rate we get 
$$\frac{\log T_n}{T_n} = \log (n \Delta_n) \frac{1}{n} n^{\frac{1}{\bar{\beta} + 1}} \frac{1}{\log(n \Delta_n)} = (\frac{1}{n})^{\frac{ \bar{\beta}}{ \bar{\beta} + 1}},$$
which is the convergence rate obtained in the intermediate regime, {\rew as $\frac{2 \bar{\beta}}{2 \bar{\beta} + 2} = \frac{ \bar{\beta}}{ \bar{\beta} + 1}$.}}\\
\\
In this section we have found the convergence rates for the estimation of the invariant density starting from the discrete observation of the process $X$. Such observations are, in this section, all taken at the same instant. One may wonder if it is possible to recover the same results when the process $X$ is observed asynchronously. The goal of next section is to answer {\rev this} question.

\section{Main results, asynchronous framework}{\label{s: asynch}}

In this section we assume $d \ge 3$ and we suppose that the components of the process $X$ are observed asynchronously, i.e. in different instants. We will see that, up to require conditions on the discretization step and on the asynchronicity of the observations, it is possible to obtain the continuous convergence rate $(\frac{\log T_n}{T_n})^{\frac{2 \bar{\beta}_3}{2 \bar{\beta}_3 + d - 2}}$. \\
%\subsection{Discrete asynchronous observations}{\label{s: as discrete}}
\\
We assume that we dispose of the discrete observations $X_{t_1^l}^l, ..., X_{t_n^l}^l$ for any $l \in \{ 1, ..., d \}$, with { \rev 
	${\modarn 0 = t_0^l < {}} t_1^l < ... < t_n^l \le T_n$ and $T_n \to \infty $ for $n \rightarrow \infty$ for any $l \in \left \{1, ... , d \right \}$.} We also define 
\begin{equation*}
    \Delta_n := \sup_{l =1, ... , d} \quad \sup_{i = 0, ... , {\rev n}} (t_{i + 1}^l - t_i^l),  \text{\rev ~ where $t^l_{n+1}=T_n$ for all $l \in \{ 1, ..., d \}$.}
\end{equation*}
We remark it would have been possible to consider a different number of points on the different directions, having in particular $n_l$ observations for the coordinate $X^l$. We have decided to take $n_1 = ... = n_d$ in order to lighten the notation. \\
Before we proceed with the statements of our results, we introduce some functions. First of all we observe that we have defined $d$ {\rev partitions} of $[0, T_n]$ and so, for any $u \in [0, T_n]$, there exist some indexes $i_1, ..., i_d$ such that $u \in [t_{i_l}^l, t_{i_l + 1}^l)$, depending on the considered direction. We introduce then the following $d$ functions $\varphi_{n,l}: [0, T_n] \rightarrow \R$ such that $\varphi_{n,l}(u) := t_{i_l}^l$, for any $l \in \{ 1, ... , d \}$. 
Thanks to these functions we can write the sums in the form of integrals. It leads us to the following estimator, which is the natural adaptation of the one in \eqref{eq: def estimator}:
\begin{align}{\label{eq: estimator asynch}}
 \hat{\pi}_{h,{T_n}}^a(x) & = \frac{1}{T_n \prod_{l = 1}^d h_l} \int_0^{T_n} \prod_{l = 1}^d K(\frac{x_l - X_{\varphi_{n, l}(u)}^l}{h_l}) du \\
 & = : \frac{1}{T_n} \int_0^{T_n} \prod_{l = 1}^d K_{h_l}(x_l - X_{\varphi_{n, l}(u)}^l) du.   \nonumber
\end{align}
{\rew In this context we introduce the class of coefficient $\tilde{\Sigma}=
	\tilde{\Sigma}(\beta, \mathcal{L},a_{\text{min}},a_0,a_1,b_0,b_1,\tilde{C},\tilde{\rho})$ as 
\begin{equation}\label{E:def_Sigma_tilde}
	\tilde{\Sigma} := \{ (a, b) \in \Sigma(\beta, \mathcal{L},a_{\text{min}},a_0,a_1,b_0,b_1,\tilde{C},\tilde{\rho}): \, a,b \mbox{ are } \mathcal{C}^3 \}.
\end{equation}
The regularity requested is needed in order to obtain uniform mixing inequalities on the class of coefficients and control on the bias of the estimator. }

{\rev We provide two different bounds for the variance, depending on the regime considered. When the discretization error is negligible we have the following result.}
\begin{proposition}
Suppose that A1-A2 hold with some constant $0<a_\text{min}\le a_0$ and $a_1>0, ~ b_0>0, b_1>0, \tilde{C}_b, \tilde{\rho}_b$ and that $d \ge 3$. {\rew Let $\beta = (\beta_1, ... , \beta_d)$, $0 < \beta_1 \le ... \le \beta_d$ and $\mathcal{L} =(\mathcal{L}_1, ... , \mathcal{L}_d)$, $\mathcal{L}_i > 0$.} Let $h^* = (h_1^*, ... , h_d^*)$ be the rate optimal choice for the bandwidth $h$ given in \eqref{eq: optimal bandwidth}. We suppose moreover that $\Delta_n \le \frac{1}{4} h_1^* h_2^* $, then there exist $c > 0$ and $T_0 > 0$ such that, for $T_n \ge T_0$,
 $$Var(\hat{\pi}^a_{h^*,n}(x)) \le \frac{c}{T_n} \frac{\sum_{j = 1}^d |\log(h^*_j)|}{\prod_{l = 3}^d h^*_l}.$$   
Moreover, the constant $c$ is uniform over the set of coefficients $(a, b) \in {\rew \tilde{\Sigma}}$.
\label{prop: bound var as}
\end{proposition}
Comparing the bound here above with the results gathered in Proposition \ref{prop: bound var discrete dge3} it appears clearly that asking the condition $\Delta_n \le \frac{1}{4} h_1^* h_2^* $ is enough both in the synchronous and asynchronous frameworks to recover the same bound on the variance as in the continuous case, which is optimal for $\beta_2 < \beta_3$. \\{ \rev 
In {\rew the} intermediate regime, instead, the following proposition holds true. 
\begin{proposition}
Suppose that A1-A2 hold with some constant $0<a_\text{min}\le a_0$ and $a_1>0, ~ b_0>0, b_1>0, \tilde{C}_b, \tilde{\rho}_b$ and that $d \ge 3$. {\rew Let $\beta = (\beta_1, ... , \beta_d)$, $0 < \beta_1 \le ... \le \beta_d$ and $\mathcal{L} =(\mathcal{L}_1, ... , \mathcal{L}_d)$, $\mathcal{L}_i > 0$.} Let $\tilde{h}^* = (\tilde{h}_1^*, ... , \tilde{h}_d^*)$ be the rate optimal choice for the bandwidth $h$ in intermediate regime, given by $\tilde{h}^*_j (n) = (\frac{1}{n})^{\frac{\bar{\beta}}{\beta_j(2 \bar{\beta} + d)}}$ for any $j \in \{1, ... , d \}$. Suppose moreover that $\Delta_n \ge (\prod_{l = 1}^d \Tilde{h}^*_l)^{\frac{2}{d}} $, then there exist $c > 0$ and $T_0 > 0$ such that, for $T_n \ge T_0$,
 $$Var(\hat{\pi}^a_{\tilde{h}^*,n}(x)) \le \frac{c}{T_n} \frac{\Delta_n}{\prod_{l = 1}^d \Tilde{h}^*_l}.$$   
Moreover, the constant $c$ is uniform over the set of coefficients $(a, b) \in {\rew \tilde{\Sigma}}$.
\label{prop: bound var as intermediate}
\end{proposition}
{\modarn If $\frac{n\Delta_n}{T_n}$ is bounded it is possible to compare the conditions on the sampling step required by Propositions \ref{prop: bound var as}--\ref{prop: bound var as intermediate}. Indeed, under the condition $\frac{n\Delta_n}{T_n}\le c$, exactly as
	we obtain \eqref{eq: 8.5}, we see that the condition on $\Delta_n$ in Proposition \ref{prop: bound var as} is equivalent to 
	$\Delta_n \underset{\sim}{<}n^{-\frac{\overline{\beta}}{2\overline{\beta}+d}(\frac{1}{\beta_1}+\frac{1}{\beta_2})}$. On  the other hand, the condition on the sampling step in Proposition \ref{prop: bound var as intermediate} is always equivalent to $\Delta_n \underset{\sim}{>}n^{-\frac{2}{2\overline{\beta}+d}}$. In the isotropic case, the two exponents are the same, and thus Proposition \ref{prop: bound var as}--\ref{prop: bound var as intermediate} covers all the cases for $\Delta_n$. In the anisotropic case, there is a possible gap corresponding to $ n^{-\frac{\overline{\beta}}{2\overline{\beta}+d}(\frac{1}{\beta_1}+\frac{1}{\beta_2})}<\Delta_n<n^{-\frac{2}{2\overline{\beta}+d}} $, where the sampling step is not small enough to recover the continuous rate and conditions required for Proposition \ref{prop: bound var as intermediate} are not satisfied either.}\\
\\
In order to obtain the convergence rate in {\modarn the asynchronous framework} we have to provide a bound on the bias term. {\rew While} in the synchronous context this object has already been deeply studied, this is no longer the case when the process is observed in different instants of time. Hence, in the asynchronous case we introduce a quantity whose goal is to measure the asynchronicity of our data:
{\modarn 
	\begin{equation} \label{E: def Delta prime}
		\Delta_n' := \sup_{t \in [0, T_n]} \sup_{i, j = 1, ... , d} |\varphi_{n, i}(t) - \varphi_{n,j}(t)|.
	\end{equation}
}
{\modarn By definition we have $\Delta_n'\le \Delta_n$, and $\Delta_n'=0$ for synchronous data.}
Then, the following proposition holds true. }
\begin{proposition}\label{prop: bias asynch}
	Suppose that A1 holds {\rew and that $d \ge 3$.} {\rew Let $\beta = (\beta_1, ... , \beta_d)$, $0 < \beta_1 \le ... \le \beta_d$ and $\mathcal{L} =(\mathcal{L}_1, ... , \mathcal{L}_d)$, $\mathcal{L}_i > 0$.} Then, there exists $c>0$ such that for all  $T_n >0$, $0<h_i<1$,
$$	\abs{\E[ \hat{\pi}_{h,{T_n}}^a(x)]- \pi(x)} \le c \sum_{i=1}^d h_i^{\beta_i} + c \sqrt{\Delta_n'}.$$	
Moreover, the constant $c$ is uniform over the set of coefficients $(a, b) \in {\rew \tilde{\Sigma}}$.
\end{proposition}
{\modarn We see that asynchronicity introduces an additional term in the control of the {\rew bias}.}
From Propositions \ref{prop: bound var as} and \ref{prop: bias asynch} {\rev the} next theorem easily follows.
\begin{theorem}
%	Suppose that A1-A2 hold {\rew and that $d \ge 3$.} {\rew Let $\beta = (\beta_1, ... , \beta_d)$, $0 < \beta_1 \le ... \le \beta_d$ and $\mathcal{L} =(\mathcal{L}_1, ... , \mathcal{L}_d)$, $\mathcal{L}_i > 0$.
{\rew Suppose that $d \ge 3$.}
{\rew Let $\beta = (\beta_1, ... , \beta_d)$, $0 < \beta_1 \le ... \le \beta_d$, $\mathcal{L} =(\mathcal{L}_1, ... , \mathcal{L}_d)$, $\mathcal{L}_i > 0$, $0<a_\text{min}\le a_0$,
	$a_1>0, ~ b_0>0, ~b_1>0, ~\tilde{C}>0, ~\tilde{\rho}>0$ and
	 $\tilde{\Sigma}$ given by \eqref{E:def_Sigma_tilde}}.		
	If {\rev $\Delta_n \le (\frac{\log T_n}{ T_n})^{ \frac{\bar{\beta}_3}{2\bar{\beta}_3 + d - 2}(\frac{1}{\beta_1} + \frac{1}{\beta_2})}$ and  $\Delta_n' \le (\frac{\log T_n}{T_n})^{\frac{2 \bar{\beta}_3}{2 \Bar{\beta}_3 + d - 2}}$}, then {\rew there exist %$\tilde{C}$, $\tilde{\rho}$, 
		$c > 0$ and $T_0 > 0$ such that, for $T_n \ge T_0$,}
$$\sup_{(a, b) \in {\rew \tilde{\Sigma}}} \mathbb{E}[|\hat{\pi}^a_{h,{T_n}}(x) - \pi (x)|^2] {\rew \le} \,
    c (\frac{\log T_n}{ T_n})^{ \frac{2\bar{\beta}_3}{2\bar{\beta}_3 + d - 2}}.$$
\label{th: conv rate asynch}
\end{theorem}
{\rev Regarding the intermediate regime, instead, Proposition \ref{prop: bound var as intermediate} and \ref{prop: bias asynch} yields the following 
\begin{theorem}
	{\rew Suppose that $d \ge 3$.}
	{\rew Let $\beta = (\beta_1, ... , \beta_d)$, $0 < \beta_1 \le ... \le \beta_d$, $\mathcal{L} =(\mathcal{L}_1, ... , \mathcal{L}_d)$, $\mathcal{L}_i > 0$, $0<a_\text{min}\le a_0$,
		$a_1>0, ~ b_0>0, ~b_1>0, ~\tilde{C}>0, ~\tilde{\rho}>0$ and
		$\tilde{\Sigma}$ given by \eqref{E:def_Sigma_tilde}}.		
	%Suppose that A1-A2 hold and {\rew that $d \ge 3$.} {\rew Let $\beta = (\beta_1, ... , \beta_d)$, $0 < \beta_1 \le ... \le \beta_d$ and $\mathcal{L} =(\mathcal{L}_1, ... , \mathcal{L}_d)$, $\mathcal{L}_i > 0$.}  
	If {\rev $\Delta_n \ge (\frac{1}{ n})^{ \frac{2}{2\bar{\beta} + d}}$ and {\rew $\Delta_n' \le (\frac{1}{n})^{ \frac{2\bar{\beta}}{2\bar{\beta} + d}}$}}, then {\rew there exist %$\tilde{C}$, $\tilde{\rho}$, 
		$c > 0$ and 
	$T_0 > 0$ such that, for $T_n \ge T_0$,}
$$\sup_{(a, b) \in {\rew \tilde{\Sigma}}} \mathbb{E}[|\hat{\pi}^a_{h,{T_n}}(x) - \pi (x)|^2] {\rew \le} \,
    c (\frac{1}{n})^{ \frac{2\bar{\beta}}{2\bar{\beta} + d}}.$$
\label{th: conv rate asynch intermediate}
\end{theorem}
Theorems \ref{th: conv rate asynch} and \ref{th: conv rate asynch intermediate} extend to the asynchronous case Theorems \ref{th: discrete d ge 3} and \ref{th: discrete d ge 3 non negl}, respectively. One can see that the bounds on the variance appearing in Proposition \ref{prop: bound var as} are stronger than the {\rew one} presented in Proposition \ref{prop: bound var as}, as they allow us to recover the optimal continuous convergence rate, while Theorem \ref{th: conv rate asynch} provides the optimal convergence rate only up {\modarn to} a logarithmic factor. In intermediate regime, instead, we obtain in the asynchronous framework the very same result as in the synchronous one. }

%\sout{The theorem here above relies on the fact that the condition on the discretization step gathered in Proposition \ref{prop: bound var as} is negligible compared to the one in Proposition \ref{prop: bias asynch}. Indeed, after having replaced the value of the optimal bandwidth as in \eqref{eq: optimal bandwidth}, we have that 
%$$h_1^* h_2^*  = (\frac{\log T_n}{ T_n})^{ \frac{\bar{\beta}_3}{2\bar{\beta}_3 + d - 2}(\frac{1}{\beta_1} + \frac{1}{\beta_2})} .$$
%As $\beta_i \ge 1$ $\forall i$, it is $\frac{1}{\beta_1} + \frac{1}{\beta_2} \le 2$, from which we derive that 
%$$(\frac{\log T_n}{ T_n})^{ \frac{\bar{\beta}_3}{2\bar{\beta}_3 + d - 2}(\frac{1}{\beta_1} + \frac{1}{\beta_2})}  > (\frac{\log T_n}{ T_n})^{ \frac{2 \bar{\beta}_3}{2\bar{\beta}_3 + d - 2}}.$$ 
%Hence, if $\Delta_n \le (\frac{\log T_n}{T_n})^{\frac{2 \bar{\beta}_3}{2 \Bar{\beta}_3 + d - 2}}$, then it is also $\Delta_n \le h_1^* h_2^*$.}} 
{ \rev It is interesting to remark that, in the case where the two components associated to the two smallest smoothnesses are continuously observed, it is possible to significantly lighten the condition on the discretization step needed to recover the bound on the variance as in Proposition \ref{prop: bound var as}.} %\subsection{What if two out of d components are continuously observed?}{\label{s: two as}}
In this case, the continuous trajectories of $X^1$ and $X^2$ are available, as well as the discrete observations $X_{t_1^l}^l, ..., X_{t_n^l}^l$ for any $l \in \{ 3, ..., d \}$, with $0\le t_1^l \le ... \le t_n^l \le T_n$. The discretization step is defined as 
{\rew \begin{equation*}
    \Delta_n := \sup_{l =3, ... , d} \quad \sup_{i = 0, ... ,n} (t_{i + 1}^l - t_i^l),  \text{\rev ~ {\rew where $t^l_{n+1}=T_n$ for all $l \in \{ 3, ..., d \}$.}}
\end{equation*}}
{\rev The kernel estimator in this context is the following:}
\begin{align}{\label{eq: estimator 2comp}}
 \bar{\pi}_{h,{T_n}}(x) & = \frac{1}{T_n \prod_{l = 1}^d h_l} \int_0^{T_n} \prod_{m = 1,2} K(\frac{x_m - X_u^m}{h_m}) \prod_{l = 3}^d K(\frac{x_l - X_{\varphi_{n, l}(u)}^l}{h_l}) du \\
 & = : \frac{1}{T_n} \int_0^{T_n} \prod_{m = 1,2} K_{h_m} (x_m - X_u^m) \prod_{l = 3}^d K_{h_l}(x_l - X_{\varphi_{n, l}(u)}^l) du.  \nonumber 
\end{align}
{\rev Then, it is possible to recover the same upper bound on the variance as in Proposition 2 of \cite{Companion}, where all the components of the process $X$ were continuously available.}
\begin{proposition}
Suppose that A1-A2 hold and that {\rew $d \ge 3$}. {\rew Let $\beta = (\beta_1, ... , \beta_d)$, $0 < \beta_1 \le ... \le \beta_d$ and $\mathcal{L} =(\mathcal{L}_1, ... , \mathcal{L}_d)$, $\mathcal{L}_i > 0$.}
Let $h^* = (h_1^*, ... , h_d^*)$ be the rate optimal choice for the bandwidth $h$ given in \eqref{eq: optimal bandwidth} and suppose that $\Delta_n \le \frac{1}{2} (\prod_{l \ge 3} h_l^*)^{\frac{2}{d - 2}} = \frac{1}{2} (\frac{\log T_n}{T_n})^{\frac{2}{2 \bar{\beta}_3 + d - 2}}$, then there exist $c > 0$ and $T_0 > 0$ such that, for $T_n \ge T_0$,
   $$Var(\bar{\pi}_{h^*,{T_n}}(x)) \le  \frac{c}{T_n} \frac{\sum_{j = 1}^d |\log(h^*_j)|}{\prod_{l = 3}^d h^*_l}.$$ 
   Moreover, the constant $c$ is uniform over the set of coefficients $(a,b) \in {\rew \tilde{\Sigma}}$. 
\label{prop: var due comp continue}
\end{proposition}
{\rev Comparing the proposition here above with Proposition 2 of \cite{Companion} one can see that, when  $k_0 = 1, 2$, not having the continuous record of the last $(d-2)$ components does not interfere in the computations of the upper bound of the variance.
The condition appearing here above is less restrictive than the one in Proposition \ref{prop: bound var discrete dge3} and Theorem \ref{th: discrete d ge 3}, as $ h_1^* h_2^* \sum_{j = 1}^d |\log h_j^*| \underset{\sim}{<} (\prod_{l \ge 3} h_l^*)^{\frac{2}{d - 2}}$. Indeed, it is equivalent to ask $(\frac{\log T_n}{T_n})^{\frac{\bar{\beta}_3}{2 \bar{\beta}_3 + d - 2}(\frac{1}{\beta_1} + \frac{1}{\beta_2})} \log T_n \underset{\sim}{<}(\frac{\log T_n}{T_n})^{\frac{2}{2 \bar{\beta}_3 + d - 2}}$, which holds true in a pure anisotropic context.} \\ %It happens as $\bar{\beta}_3 (\frac{1}{\beta_1} + \frac{1}{\beta_2}) \frac{1}{2} > 1$.} \\
\\
{\rev However, it is easy to see that Proposition \ref{prop: bias asynch} still holds true when $\varphi_{n,l}(t) = t$ for $l=1, 2$ and so with $\bar{\pi}_{h,{T_n}}(x)$ instead of $\hat{\pi}^a_{h,{T_n}}(x)$. Then, the condition on $\Delta_n'$ turns out being in this case a condition on $\Delta_n$, as two {\rew components} are observed continuously. This constraint is stronger than both the conditions in Propositions \ref{prop: bound var as} and \ref{prop: var due comp continue}. It implies that one can propose an estimator based on the continuous observations of two components in order to improve the bound on the variance, but this implies a big deterioration on the condition on the bias. In particular, the resulting estimation is not better than the one based on the discrete asynchronous observation of all the components.  \\
One may observe that, even if in the paper {\rew it} is always assumed that the target stationary distribution $\pi$ is such that the regularities are ordered, to build our estimators as in \eqref{eq: def discrete estimator} and \eqref{eq: estimator asynch} it is not necessary to know in advance which coordinates are less regular. It is no longer the case for the estimator proposed in \eqref{eq: estimator 2comp}, where the two continuously observed coordinates are the less regular ones. This is the biggest limitation of Proposition \ref{prop: var due comp continue}, which is only a theoretical result. Its purpose is to remark that having finer observations in some coordinates yields a worse condition on the discretization step, which is surprising. This suggests it would be better to synchronise the continuous observations of the first two components with the closest discrete observation, to decrease the asynchronicity.} \\
The proofs of all the theorems and propositions stated in this section can be found in Section \ref{S: proof asynch}, but for Proposition \ref{prop: var due comp continue}, whose proof can be found in the appendix.

\section{Proof main results, synchronous framework}{\label{S: proofs synch}}

This section is devoted to the proof of our main results in the case where all the components are observed at the same time. {\rev In the sequel the constant $c$ may change from line to line.}

\subsection{Proof of Proposition \ref{prop: bound var discrete dge3}}
\begin{proof}
The proof of Proposition \ref{prop: bound var discrete dge3} heavily relies on the proof of the upper bound on the variance of \eqref{eq: def estimator}, in the continuous case. Intuitively, the integrals in Proposition 2 of \cite{Companion} will be now replaced by sums, that we will split in order to use some different bounds on each of them. The main tools are the exponential ergodicity of the process as gathered in Proposition 1 of \cite{Companion} and a bound on the transition density as in Proposition 5.1 of \cite{Gobet LAMN}.
From the definition of our estimator $\hat{\pi}_{h,n}$, using also the fact that we are considering a uniform discretization step, it follows
\begin{align*}
    Var(\hat{\pi}_{h,n}(x)) & = Var(\frac{1}{n \Delta_n} \sum_{j = 0}^{n - 1} \K_h(x - X_{t_j}) \Delta_n) \\
    & = \frac{\Delta_n^2}{T_n^2} \sum_{j = 0}^{n - 1} (n - j) \, Cov(\K_h(x-X_0), \K_h(x - X_{t_j}) ) \\
    & = : \frac{\Delta_n^2}{T_n^2} (\sum_{j = 0}^{j_{\delta_1}} + \sum_{j = j_{\delta_1} + 1}^{j_{\delta_2}} + \sum_{j = j_{\delta_2} +1 }^{j_D} + \sum_{j = j_D + 1}^{n - 1} )\,(n - j) \,  k(t_j) \\
    & =: I_1 + I_2 + I_3 + I_4,
\end{align*}
having introduced $0 \le j_{\delta_1} \le j_{\delta_2} \le j_D \le n-1$ and set $\delta_1 = \Delta_n j_{\delta_1}$, $\delta_2 = \Delta_n j_{\delta_2}$, $D = \Delta_n j_{D}$ and 
$$k(t) := Cov(\K_h(x-X_0), \K_h(x - X_{t}) ).$$
Recall also that $t_j = \Delta_n j$. 
The quantities $j_{\delta_1}, \, j_{\delta_2}, \, j_D$ (and consequently $\delta_1$, $\delta_2$ and $D$) will be chosen later, in order to get an upper bound on the variance as sharp as possible.
We will provide some bounds for $I_1$, and $I_4$ which do not depend on $k_0$, while we will bound differently $I_2$ and $I_3$ depending on whether or not $k_0$ is larger than $3$. For $j$ small we use Cauchy-Schwarz inequality, the stationarity of the process, the boundedness of $\pi$ and the definition of the kernel function to obtain
\begin{equation}
|k(t_j)| \le Var(\mathbb{K}_h(x -X_0))^\frac{1}{2}Var(\mathbb{K}_h(x -X_{t_j}))^\frac{1}{2} \le \int_{\mathbb{R}^d} (\mathbb{K}_h (x -y))^2 \pi(y) dy \le \frac{c}{\prod_{l = 1}^d h_l}.
\label{eq: bound k 5.5}
\end{equation}
It follows
\begin{equation}
|I_1| \le \frac{\Delta_n^2 \, n}{T_n^2} \sum_{j = 0}^{j_{\delta_1}} \frac{c}{\prod_{l = 1}^d h_l} = c \frac{\Delta_n^2 \, n}{T_n^2}  \frac{1}{\prod_{l = 1}^d h_l}(j_{\delta_1} + 1).
\label{eq: bound I1}
\end{equation}
For $j \in [j_{\delta_1} + 1, j_{\delta_2}]$ we act differently depending on $k_0$ as done for $s\in [\delta_1, \delta_2)$ in Proposition 2 of \cite{Companion}. For $k_0 = 1$ and $\beta_2 < \beta_3$ or $k_0 = 2$ it provides (see Equation (15) in \cite{Companion})
$$|k(s)| \le \frac{c}{\prod_{j \ge 3} h_j} \frac{1}{s}$$
that, for $s= t_j$, becomes
$$|k(t_j)| \le \frac{c}{\prod_{l \ge 3} h_l} \frac{1}{t_j}.$$
It yields
$$|I_2| \le \frac{\Delta_n^2 \, n}{T_n^2} \sum_{j = j_{\delta_1} + 1}^{j_{\delta_2}} \frac{c}{\prod_{l \ge 3} h_l} \frac{1}{t_j} = c \frac{\Delta_n}{T_n} \frac{1}{\prod_{l \ge 3} h_l}  \sum_{j = j_{\delta_1} + 1}^{j_{\delta_2}} \frac{1}{t_j}.$$
We recall that $t_j$ can be seen as $\Delta_n \, j$. Therefore, 
\begin{align*}
\sum_{j = j_{\delta_1} + 1}^{j_{\delta_2}} \frac{1}{t_j} \le \frac{c}{\Delta_n} \log(\frac{j_{\delta_2}}{j_{\delta_1}}) = \frac{c}{\Delta_n} \log(\frac{\Delta_n j_{\delta_2}}{\Delta_n j_{\delta_1}}) = \frac{c}{\Delta_n} \log(\frac{\delta_2}{\delta_1}).
\end{align*}
It follows 
\begin{equation}
|I_2| \le  c \frac{1}{T_n} \frac{1}{\prod_{l \ge 3} h_l} \log(\frac{\delta_2}{\delta_1}).
\label{eq: bound I2}
\end{equation}
When $k_0 \ge 3$ instead, acting as to get (17) in \cite{Companion} and taking $s = t_j$, we obtain 
$$|k(t_j)| \le \frac{c}{\prod_{l \ge k_0 + 1} h_l} t_j^{- \frac{k_0}{2}}.$$
Therefore, 
\begin{align}{\label{eq: bound I2 k0 large}}
|I_2| & \le \frac{\Delta_n^2 \, n}{T_n^2} \sum_{j = j_{\delta_1} + 1}^{j_{\delta_2}} \frac{c}{\prod_{l \ge k_0 + 1} h_l} t_j^{- \frac{k_0}{2}} \nonumber \\
& = c \frac{\Delta_n}{T_n} \frac{1}{\prod_{l \ge k_0 + 1} h_l}  \sum_{j = j_{\delta_1} + 1}^{j_{\delta_2}} \Delta_n^{- \frac{k_0}{2}} j^{- \frac{k_0}{2}} \nonumber \\
& \le c \frac{\Delta_n^{1 - \frac{k_0}{2}}}{T_n} \frac{1}{\prod_{l \ge k_0 + 1} h_l} j_{\delta_1}^{1 - \frac{k_0}{2}} \nonumber \\
&= c \frac{\delta_1^{1 - \frac{k_0}{2}}}{T_n}\frac{1}{\prod_{l \ge k_0 + 1} h_l} 
\end{align}
where we have used that, as $k_0 \ge 3$, $1 - \frac{k_0}{2}$ is negative. \\
To conclude the analysis of $I_2$ we assume that $k_0 = 1$ and $\beta_2 = \beta_3$. In this case the estimation here above still holds but, as $1 - \frac{k_0}{2} = \frac{1}{2}$ is now positive, it provides 
\begin{equation}
|I_2| \le c \frac{\delta_2^{\frac{1}{2}}}{T_n} \frac{1}{\prod_{l \ge 2} h_l}.
\label{eq: I2 beta2= beta3 39.5}
\end{equation}
We now deal with $I_3$. With the same bound on the covariance as in (20) of Proposition 2 in \cite{Companion} we get in any case, but for $k_0 = 1$ and $\beta_2 = \beta_3$,
$$| k(s)| \le c (s^{- \frac{d}{2}} + 1).$$
Therefore, taking $s = t_j$,
\begin{align}{\label{eq: bound I3}}
|I_3| & \le \frac{\Delta_n^2\, n}{T_n^2} \sum_{j = j_{\delta_2} +1 }^{j_D} c (t_j^{- \frac{d}{2}} + 1) \\
& \le c \frac{\Delta_n}{T_n}(  \sum_{t_j \le 1, \, j = j_{\delta_2} +1 }^{j_D} \Delta_n^{- \frac{d}{2}}  j^{- \frac{d}{2}} +  \sum_{t_j > 1, \, j = j_{\delta_2} +1 }^{j_D} 1) \nonumber \\
& \le c \frac{\Delta_n}{T_n} (\Delta_n^{- \frac{d}{2}}j_{\delta_2 + 1}^{1 - \frac{d}{2}}  + j_D) \nonumber \\
& = \frac{c}{T_n} (\delta_2^{1 - \frac{d}{2}} + D). \nonumber
\end{align}
For $k_0 = 1$ and $\beta_2 = \beta_3$, instead, (22) of \cite{Companion} provides 
$$|k(s)| \le c (s^{- \frac{3}{2}} \frac{1}{\prod_{l \ge 4} h_l} + 1).$$
It follows
$$|I_3| \le \frac{c \Delta_n^2 n}{T_n^2} \frac{1}{\prod_{l \ge 4} h_l} \sum_{j = j_{\delta_2} +1 }^{j_D} (t_j^{- \frac{3}{2}} + 1). $$
Acting as above we obtain
\begin{equation}
|I_3| \le \frac{c}{T_n} (\frac{1}{\prod_{l \ge 4} h_l} \frac{1}{\delta_2^{\frac{1}{2}}} + D).
\label{eq: I3 beta2 = beta3 40.5}
\end{equation}
To conclude, we need to evaluate the case where $j \in [j_D + 1, n-1]$. In this interval we use the exponential ergodicity of the process, as in Proposition 1 of \cite{Companion}. It follows 
$$|k(t_j)| \le c \left \| \mathbb{K}_h (x - \cdot ) \right \|_\infty^2 e^{- \rho t_j} \le \frac{c}{(\prod_{l = 1}^d h_l)^2}e^{- \rho t_j},$$
for $c$ and $\rho$ positive constant uniform over the set of coefficients $(a, b) \in \Sigma$.
It implies 
\begin{align}{\label{eq: bound I4}}
|I_4| & \le c \frac{\Delta_n}{T_n} \frac{1}{(\prod_{l = 1}^d h_l)^2} \sum_{j = j_D + 1}^{n - 1} e^{- \rho \Delta_n j} \\
& \le c \frac{\Delta_n}{T_n} \frac{1}{(\prod_{l = 1}^d h_l)^2}e^{- \rho \Delta_n (j_D + 1)} \nonumber \\
& \le c \frac{\Delta_n}{T_n} \frac{1}{(\prod_{l = 1}^d h_l)^2}e^{- \rho D} \nonumber. 
\end{align}
From \eqref{eq: bound I1}, \eqref{eq: bound I2}, \eqref{eq: bound I3} and \eqref{eq: bound I4} we obtain the following bound for the case $k_0 = 1$ and $\beta_2< \beta_3 $ or $k_0 = 2$: 
\begin{align*}
Var(\hat{\pi}_{h,n}(x)) & \le \frac{c}{T_n} \big[ \frac{1}{\prod_{l =1}^d h_l}(\delta_1 + \Delta_n) + \frac{1}{\prod_{l \ge 3} h_l} \log (\frac{\delta_2}{\delta_1}) + \delta_2^{1 - \frac{d}{2}} + D + \frac{1}{(\prod_{l = 1}^d h_l)^2}e^{- \rho D} \big].
\end{align*}
It is easy to see that, except for the term $\frac{c}{T_n} \frac{1}{\prod_{l =1}^d h_l} \Delta_n$, the bound is the same as in the continuous case (see (22) in \cite{Companion}). Hence, also the optimal choice for the parameters $\delta_1$, $\delta_2$ and $D$ should be the same as in the continuous case, for which $\delta_1 = h_1 h_2$, $\delta_2 := (\prod_{j \ge 3} h_j)^{\frac{2}{d-2}}$ and $D:= [\max (- \frac{2}{\rho} \log (\prod_{j = 1}^d h_j), 1) \land T]$.
Recalling that $j_{\delta_1}$, $j_{\delta_2}$ and $j_D$ have to be some integers, we can not propose exactly the same choice as above but we can take 
$$j_{\delta_1} := \lfloor \frac{h_1 h_2}{\Delta_n} \rfloor, \qquad j_{\delta_2} := \lfloor \frac{(\prod_{j \ge 3} h_j)^{\frac{2}{d-2}}}{\Delta_n} \rfloor, \qquad j_{D} := \lfloor \frac{[\max (- \frac{2}{\rho} \log (\prod_{j = 1}^d h_j), 1) \land T]}{\Delta_n} \rfloor.$$
It yields 
$$Var(\hat{\pi}_{h,n}(x)) \le \frac{c}{T_n}\frac{\sum_{l = 1}^d |\log(h_l)|}{\prod_{l \ge 3} h_l} + \frac{c}{T_n} \frac{1}{\prod_{l =1}^d h_l} \Delta_n, $$
as we wanted. \\
When $k_0 \ge 3$, instead, we replace the bound gathered in \eqref{eq: bound I2} with the one in \eqref{eq: bound I2 k0 large}. It follows
\begin{align*}
Var(\hat{\pi}_{h,n}(x)) & \le \frac{c}{T_n} \big[ \frac{1}{\prod_{l =1}^d h_l}(\delta_1 + \Delta_n) + \frac{1}{\prod_{l \ge k_0 + 1} h_l} \delta_1^{1 - \frac{k_0}{2}}  + \delta_2^{1 - \frac{d}{2}} + D + \frac{1}{(\prod_{l = 1}^d h_l)^2}e^{- \rho D} \big].
\end{align*}
Again, every term but $\frac{c}{T_n} \frac{1}{\prod_{l =1}^d h_l} \Delta_n$ was already present in the proof of Proposition 2 of \cite{Companion} (see (23)) and so the best choice would be to take the parameters as before. With this purpose in mind we choose
$$j_{\delta_1} := \lfloor \frac{(\prod_{l = 1}^{k_0} h_l)^{\frac{2}{k_0}}}{\Delta_n} \rfloor, \qquad j_{\delta_2} := \lfloor \frac{1}{\Delta_n} \rfloor, \qquad j_{D} := \lfloor \frac{[\max (- \frac{2}{\rho} \log (\prod_{j = 1}^d h_j), 1) \land T]}{\Delta_n} \rfloor.$$
It follows
$$Var(\hat{\pi}_{h,n}(x)) \le \frac{c}{T_n}\frac{1}{(\prod_{l = 1}^{k_0} h_l)^{1 - \frac{2}{k_0}}(\prod_{l \ge k_0 +1} h_l)} + \frac{c}{T_n} \frac{1}{\prod_{l =1}^d h_l} \Delta_n. $$
We are left to study the case where $k_0 = 1$ and $\beta_2 = \beta_3$. Here, from \eqref{eq: bound I1}, \eqref{eq: I2 beta2= beta3 39.5}, \eqref{eq: I3 beta2 = beta3 40.5} and \eqref{eq: bound I4} we obtain 
\begin{align*}
Var(\hat{\pi}_{h,n}(x)) & \le \frac{c}{T_n} \big[ \frac{1}{\prod_{l =1}^d h_l}(\delta_1 + \Delta_n) + \frac{\delta_2^{\frac{1}{2}}}{\prod_{l \ge 2} h_l}  + \frac{1}{\prod_{l \ge 4} h_l \, \delta_2^{\frac{1}{2}}} + D + \frac{1}{(\prod_{l = 1}^d h_l)^2}e^{- \rho D} \big].
\end{align*}
We take 
$$j_{\delta_1} := 1, \qquad j_{\delta_2} := \lfloor \frac{h_2 h_3}{\Delta_n} \rfloor, \qquad j_{D} := \lfloor \frac{[\max (- \frac{2}{\rho} \log (\prod_{j = 1}^d h_j), 1) \land T]}{\Delta_n} \rfloor$$
to get
$$Var(\hat{\pi}_{h,n}(x)) \le \frac{c}{T_n}\frac{1}{\sqrt{h_2 h_3}\prod_{l \ge 4} h_l} + \frac{c}{T_n} \frac{1}{\prod_{l =1}^d h_l} \Delta_n. $$
All the constant are uniform over the set of coefficients $(a,b) \in \Sigma$.
The proof of Proposition \ref{prop: bound var discrete dge3} is then complete.
\end{proof}

\subsection{Proof of Theorem \ref{th: discrete d ge 3}}
\begin{proof}
{\rev We split the proof according to the choice of $k_0$.} \\
$\bullet$ For $k_0 = 1$ and $\beta_2 < \beta_3$ or $k_0 = 2$, when $\Delta_n \le h_1^* h_2^* \sum_{j = 1}^d |\log h_j^*| = (\frac{\log T_n}{T_n})^{\frac{\bar{\beta}_3}{2 \bar{\beta}_3 + d - 2}(\frac{1}{\beta_1} + \frac{1}{\beta_2})} \log T_n$ the result is a straightforward consequence of the bias variance decomposition and of the bound on the variance gathered in Proposition \ref{prop: bound var discrete dge3}. The rate optimal choice of the bandwidth $h^*$ as in \eqref{eq: optimal bandwidth} provides 
\begin{align*}
 Var(\hat{\pi}_{h,n}(x)) & \le \frac{c}{T_n}\frac{\sum_{l = 1}^d |\log(h^*_l)|}{\prod_{l \ge 3} h^*_l} + \frac{c}{T_n} \frac{1}{\prod_{l =1}^d h_l^*} \Delta_n  \\
 & \le \frac{c}{T_n}\frac{\sum_{l = 1}^d |\log(h^*_l)|}{\prod_{l \ge 3} h^*_l}. 
\end{align*}
Hence, 
\begin{align*}
\sup_{(a,b) \in \Sigma} \mathbb{E}[|\hat{\pi}_{h,n}(x) - \pi (x)|^2] & \le c \sum_{j = 1}^d h_j^{* \, 2 \beta_j} + \frac{c}{T_n}\frac{\sum_{l = 1}^d |\log(h^*_l)|}{\prod_{l \ge 3} h^*_l} \\
& = (\frac{\log T_n}{T_n})^{\frac{\bar{\beta}_3}{(2 \bar{\beta}_3 + d - 2)}}
\end{align*}
 $\bullet$ For $k_0 \ge 3$, when $\Delta_n  \le (h_1^* h_2^*) = (h_1^*)^2 = (\frac{1}{T_n})^{\frac{2 \bar{\beta}_3}{\beta_1 (2 \bar{\beta}_3 + d - 2)}}$, we get the continuous convergence rate by the bias-variance decomposition and the second point of Proposition \ref{prop: bound var discrete dge3}, choosing $h_1^* = ... = h^*_{k_0}$ and $h_l^*(T_n) : (\frac{1}{T_n})^{\frac{\bar{\beta}_3}{\beta_l (2 \bar{\beta}_3 + d - 2)}}$ for any $l \in \{ 1, ... , d \}$.   \\
 The reasoning is the same for $k_0 = 1$ and $\beta_2 = \beta_3$, recalling that we no longer have $h_1^* = h_2^*$ (as $\beta_1 < \beta_2$) and so the condition on the discretization step becomes
 $\Delta_n  \le h_1^* h_2^* =: (\frac{1}{T_n})^{\frac{ \bar{\beta}_3}{2 \bar{\beta}_3 + d - 2}(\frac{1}{\beta_1} + \frac{1}{\beta_2})}$.
\end{proof}

\subsection{Proof of Theorem \ref{th: discrete d ge 3 non negl}}
\begin{proof}
Even if the final result is the same for any $k_0$, the proof of Theorem \ref{th: discrete d ge 3 non negl} is substantially different depending on the ordering of the smoothness and so on the value of $k_0$.  \\
$\bullet$ We start assuming $k_0 =1$ and $\beta_2 < \beta_3$ or $k_0 = 2$. When $\Delta_n > (\frac{\log T_n}{T_n})^{\frac{\bar{\beta}_3}{2 \bar{\beta}_3 + d - 2}(\frac{1}{\beta_1} + \frac{1}{\beta_2})} \log T_n$, the form of the rate optimal bandwidth is no longer as in the continuous case. We observe that $T_n = n \Delta_n$ and
$$\frac{1}{\beta_1} + \frac{1}{\beta_2} = \frac{d}{\bar{\beta}}- \frac{d-2}{\bar{\beta_3}},$$
and so the condition here above is equivalent to ask {\rev $\Delta_n^{1 + \alpha} > (\log T_n)^{1 + \alpha} (\frac{1}{n})^\alpha$. It holds true if and only if}
\begin{equation}
\Delta_n > (\log T_n) (\frac{1}{n})^{\frac{\alpha}{1 + \alpha}},
\label{eq: con Delta}
\end{equation}
with 
\begin{equation}
\alpha:=\frac{\bar{\beta}_3}{2 \bar{\beta}_3 + d - 2}(\frac{1}{\beta_1} + \frac{1}{\beta_2}) = \frac{\bar{\beta}_3}{2 \bar{\beta}_3 + d - 2}(\frac{d}{\bar{\beta}}- \frac{d-2}{\bar{\beta_3}}) = \frac{\bar{\beta}_3 d - (d-2) \bar{\beta}}{\bar{\beta} (2 \bar{\beta}_3 + d - 2)}.
\label{eq: def alpha}
\end{equation}
From the bias variance decomposition together with Proposition \ref{prop: bound var discrete dge3} we now obtain 
$$\mathbb{E}[|\hat{\pi}_{h,n}(x) - \pi (x)|^2] \le c \sum_{j = 1}^d h_j^{2 \beta_j} + \frac{c}{T_n}\frac{\sum_{l = 1}^d |\log(h_l)|}{\prod_{l \ge 3} h_l} + \frac{c}{T_n} \frac{1}{\prod_{l =1}^d h_l} \Delta_n.$$
We look for the rate optimal choice of the bandwidth by choosing $a_1$, ... , $a_d$ such that $h_l = (\frac{1}{n})^{a_l}$. We get the following bound 
\begin{align}{\label{eq: mse k0 large}}
\mathbb{E}[|\hat{\pi}_{h,n}(x) - \pi (x)|^2] & \le c \sum_{j = 1}^d (\frac{1}{n})^{2 \beta_j a_j} + \frac{c}{n \Delta_n} \frac{\sum_{l = 1}^d a_l \log n}{\prod_{ l \ge 3} (\frac{1}{n})^{a_l}} + \frac{c}{n} \frac{1}{\prod_{l = 1}^d (\frac{1}{n})^{a_l}} \\
& \le c \sum_{j = 1}^d (\frac{1}{n})^{2 \beta_j a_j} + \frac{c}{n} \frac{1}{\log (n \Delta_n) (\frac{1}{n})^{\frac{\alpha}{1 + \alpha}}} \frac{\log n}{\prod_{l \ge 3} (\frac{1}{n})^{a_l}} + \frac{c}{n} \frac{1}{\prod_{l = 1}^d (\frac{1}{n})^{a_l}}, \nonumber
\end{align}
having used the condition on $\Delta_n$ gathered in \eqref{eq: con Delta}. We recall we have assumed $n \Delta_n \rightarrow \infty$ for $n \rightarrow \infty$ and that $\Delta_n > n^{-k}$ for some $k\in (0,1)$. It follows that $\log(n \Delta_n) < (1 - k) \log n$, which implies $\frac{\log n}{\log(n \Delta_n)} \le c$. Then, we observe that the balance between the three terms in \eqref{eq: mse k0 large} is achieved for $h_l(n) = (\frac{1}{n})^{\frac{\bar{\beta}}{\beta_j(2 \bar{\beta} + d)}}$. In this way $\prod_{l \ge 3} h_l = (\frac{1}{n})^{\frac{\bar{\beta}}{2 \bar{\beta} + d} \frac{d-2}{\bar{\beta_3}}}$, which implies in particular that the second term in the right hand side of \eqref{eq: mse k0 large} is upper bounded by 
\begin{align}{\label{eq: final n mse}}
&(\frac{1}{n})^{\frac{1}{1 + \alpha} - \frac{\bar{\beta}}{2 \bar{\beta} + d} \frac{d-2}{\bar{\beta_3}}} \nonumber \\
& = (\frac{1}{n})^{\frac{\bar{\beta}(2 \bar{\beta_3} + d - 2)}{\bar{\beta_3}(2 \bar{\beta} + d)} - \frac{\bar{\beta}(d-2)}{\bar{\beta_3}(2 \bar{\beta} + d)}} \\
& = (\frac{1}{n})^{\frac{2 \bar{\beta}}{2 \bar{\beta} + d}}, \nonumber
\end{align}
 which is clearly the size of the other terms as well, after having replaced the rate optimal choice for $h_l (n)$. \\
 \\
$\bullet$ We now consider the case where $\Delta_n > (\frac{1}{T_n})^{\frac{\bar{\beta}_3}{2 \bar{\beta}_3 + d - 2}(\frac{1}{\beta_1} + \frac{1}{\beta_2})}$ and $\beta_2 = \beta_3$. We start assuming that $k_0 = 1$. We can write
$$\Delta_n > (\frac{1}{n})^{\frac{\alpha}{1 + \alpha}},$$
with $\alpha$ as in \eqref{eq: def alpha}. From the bound on the variance gathered in Proposition \ref{prop: bound var discrete dge3} and the bias-variance decomposition easily follows the wanted result, acting as above. \\
\\
When $k_0 \ge 3$ it is $\beta_1 = \beta_2$ and so $\Delta_n >  (\frac{1}{T_n})^{\frac{\bar{\beta}_3}{2 \bar{\beta}_3 + d - 2}(\frac{1}{\beta_1} + \frac{1}{\beta_2})} = (\frac{1}{T_n})^{\frac{2 \bar{\beta}_3}{\beta_1 (2 \bar{\beta}_3 + d - 2)}}$. As $T_n = n \Delta_n$, the previous condition is equivalent to ask $\Delta_n > (\frac{1}{n})^{\frac{\alpha}{1+ \alpha}}$, with 
 $$\alpha = \frac{2 \bar{\beta}_3}{\beta_1 (2 \bar{\beta}_3 + d - 2)}.$$
 We underline that $\alpha$ is exactly the same as in \eqref{eq: def alpha}, having now $\beta_1 = \beta_2$. The previous remark, together with the second point of Proposition \ref{prop: bound var discrete dge3} leads to the following bound for the mean squared error, for $k_0 \ge 3$: 
 \begin{align*}
\mathbb{E}[|\hat{\pi}_{h,n}(x) - \pi (x)|^2] & \le  c \sum_{j = 1}^d h_j^{2 \beta_j} + \frac{c}{T_n} \frac{1}{(\prod_{l = 1}^{k_0}h_l)^{1 - \frac{2}{k_0}} \prod_{l \ge k_0 + 1}h_l} + \frac{c}{T_n} \frac{\Delta_n}{\prod_{l = 1}^d h_l} \\
& \le c \sum_{j = 1}^d h_j^{2 \beta_j} + \frac{c}{n (\frac{1}{n})^{\frac{\alpha}{1 + \alpha}}} \frac{1}{(\prod_{l = 1}^{k_0}h_l)^{1 - \frac{2}{k_0}} \prod_{l \ge k_0 + 1}h_l} + \frac{c}{n} \frac{1}{\prod_{l = 1}^d h_l}. 
\end{align*}
As before, we choose the rate optimal bandwidth as $h_l(n) := (\frac{1}{n})^{\frac{\bar{\beta}}{\beta_l (2 \bar{\beta} + d)}}$. As $\beta_1 = ... = \beta_{k_0}$ it follows in particular that $h_1(n) = ... = h_{k_0} (n)$. Replacing the value of $h_l(n)$ in the bound of the mean squared error we get
$$\mathbb{E}[|\hat{\pi}_{h,n}(x) - \pi (x)|^2]  \le  c (\frac{1}{n})^{\frac{2 \bar{\beta}}{2 \bar{\beta} + d}} + c (\frac{1}{n})^{1 - \frac{\alpha}{1 + \alpha}} n^{\frac{ \bar{\beta}}{2 \bar{\beta} + d}(\frac{k_0 - 2}{\beta_1})} n^{\frac{ \bar{\beta}}{2 \bar{\beta} + d} (\frac{d- k_0}{\bar{\beta}_k})} + c (\frac{1}{n})^{1 - \frac{\bar{\beta}}{2 \bar{\beta} + d}(\frac{d}{\bar{\beta}})},$$
recalling that $\bar{\beta}_k$ is the mean smoothness over $\beta_{k_0 + 1}$, ... , $\beta_d$ and it is such that $\frac{1}{\bar{\beta}_k} = \frac{1}{d-k_0} \sum_{l \ge k_0 + 1} \frac{1}{\beta_l}$. 
We remark that 
$$\frac{k_0 - 2}{\beta_1} + \frac{d - k_0}{ \bar{\beta}_k} = \frac{d-2}{\bar{\beta}_3}.$$
Then, using also \eqref{eq: final n mse}, we have that the exponent of $\frac{1}{n}$ in the second term here above is 
$$\frac{1}{1 + \alpha} - \frac{ \bar{\beta}}{2 \bar{\beta} + d} (\frac{d-2}{\bar{\beta}_3}) = \frac{ 2 \bar{\beta}}{2 \bar{\beta} + d}.$$
Remarking that the constant $c$ does not depend on $(a,b) \in \Sigma$, it follows
$$\sup_{(a,b) \in \Sigma} \mathbb{E}[|\hat{\pi}_{h,n}(x) - \pi (x)|^2] \le c (\frac{1}{n})^{\frac{ 2\bar{\beta}}{2 \bar{\beta} + d}},$$
as we wanted.
\end{proof}

\subsection{Proof of Proposition \ref{prop: bound var discrete d=2}}
\begin{proof}
 The proof follows the procedure to bound the variance for $d=2$ proposed in the continuous case (see Theorem 2 of \cite{Companion}). We split the sum in three terms: 
\begin{align*}
	%{\label{eq: def I tilde}}
    Var(\hat{\pi}_{h,n}(x)) & = Var(\frac{1}{n \Delta_n} \sum_{j = 0}^{n - 1} \K_h(x - X_{t_j}) \Delta_n) \nonumber \\
    & = \frac{\Delta_n^2}{T_n^2} (\sum_{j = 0}^{j_{\delta}} + \sum_{j = j_{\delta} +1 }^{j_D} + \sum_{j = j_D + 1}^{n - 1} )\,(n - j) \,  k(t_j) \\
    & =: \tilde{I}_1 + \tilde{I}_2 + \tilde{I}_3. \nonumber
\end{align*}
{\rev Repeating the arguments given for the bounds on $I_1$ and $I_4$ (defined in the proof of Proposition \ref{prop: bound var discrete dge3}) we find the following bounds for $\tilde{I}_1$ and $\tilde{I}_3$, respectively.}
\begin{equation*}
|\tilde{I}_1| \le \frac{c}{T_n} \frac{1}{h_1 h_2}(\delta + \Delta_n),
%\label{eq: bound Itilde1}
\end{equation*}
\begin{equation*}
|\tilde{I}_3|\le \frac{c}{T_n} \frac{1}{(h_1 h_2)^2}e^{- \rho D}.
%\label{eq: bound Itilde3}
\end{equation*}
Regarding $\tilde{I}_2$, we act here as we did on $I_3$ in Proposition \ref{prop: bound var discrete dge3}. Recalling that here $d=2$ we get 
\begin{equation*}
|k(t_j)| \le c(t_j^{- \frac{d}{2}} + 1) \le  c( \frac{1}{t_j} + 1).
%\label{eq: bound cov low}
\end{equation*}
We need to consider separately what happens when $t_j$ is larger or smaller than $1$.
\begin{align*}
|\tilde{I}_2| & \le c \frac{\Delta_n}{T_n} \sum_{j = j_{\delta} +1 }^{j_D}(\frac{1}{t_j} + 1) \\
& \le c \frac{\Delta_n}{T_n} (\sum_{t_j \le 1, \, j = j_{\delta} +1 }^{j_D} \frac{1}{t_j} + \sum_{t_j > 1, \, j = j_{\delta} +1 }^{j_D}1) \\
& \le \frac{c}{T_n} (|\log D| + |\log \delta| + D).
\end{align*}
Putting all the pieces together, one can see that the choice $j_{\delta} := \lfloor \frac{h_1 h_2}{\Delta_n} \rfloor$ and $j_D:= \lfloor \frac{[\max (- \frac{2}{\rho} \log (h_1 h_2), 1) \land T]}{\Delta_n} \rfloor $ leads to the wanted result.
\end{proof}

\subsection{Proof of Theorem \ref{th: discrete d =2}}
\begin{proof}
The scheme we follow to prove Theorem \ref{th: discrete d =2} is the one provided in the proof of Theorem \ref{th: discrete d ge 3}. We start considering the case where 
$$\Delta_n \le h_1^* h_2^* \sum_{j = 1}^2 |\log h_j^*| = (\frac{\log T_n}{T_n})^{(\frac{1}{2 \beta_1} + \frac{1}{2 \beta_2})}\log T_n = (\frac{\log T_n}{T_n})^{\frac{1}{\bar{\beta}}} \log T_n, $$
where we have used that, from the proof of Theorem 2 of \cite{Companion}, $h_l^*(T_n) = (\frac{\log T_n}{T_n})^{a_l}$ with $a_l \ge \frac{1}{2 \beta_l}$ and that, for $d=2$, $\frac{1}{2}(\frac{1}{\beta_1} + \frac{1}{\beta_2}) = \frac{1}{\bar{\beta}}$. \\
From the bias variance decomposition together with Proposition \ref{prop: bound var discrete d=2}, taking the rate optimal choice $h_l^*(T_n)$ as above directly follows
$$\mathbb{E}[|\hat{\pi}_{h,n}(x) - \pi (x)|^2] \le c \frac{\log T_n}{T_n} + c \frac{\log(\log T_n)}{T_n} = c \frac{\log T_n}{T_n}. $$
If $\Delta_n > (\frac{\log T_n}{T_n})^{\frac{1}{\bar{\beta}}} \log T_n$, instead, it is also $\Delta_n > (\frac{1}{n})^{\frac{1}{\bar{\beta} + 1}}(\log (n \Delta_n))$. Using  the bias variance decomposition and Proposition \ref{prop: bound var discrete d=2} we obtain 
\begin{align*}
\mathbb{E}[|\hat{\pi}_{h,n}(x) - \pi (x)|^2] & \le c(h_1^{2 \beta_1} + h_2^{2 \beta_2}) + \frac{c}{T_n} \sum_{j = 1}^2 |\log h_j| + \frac{c}{T_n} \frac{\Delta_n}{h_1 h_2} \\
& \le c(h_1^{2 \beta_1} + h_2^{2 \beta_2}) + c (\frac{1}{n})^{1 - \frac{1}{\bar{\beta} + 1}} \frac{\sum_{j = 1}^2 |\log h_j|}{\log (n \Delta_n)} + \frac{c}{n h_1 h_2}.
\end{align*}
We choose the rate optimal bandwidth $h_l(n) := (\frac{1}{n})^{\frac{\bar{\beta}}{\beta_l (2 \bar{\beta} + d)}}$, for $l = 1, 2$. We have already discussed the behaviour of $n \Delta_n$, saying in particular that $\frac{\log n}{\log n \Delta_n} \le c$ in the proof of Theorem \ref{th: discrete d ge 3 non negl}. It yields 
\begin{align*}
\mathbb{E}[|\hat{\pi}_{h,n}(x) - \pi (x)|^2] & \le c(\frac{1}{n})^{\frac{2 \bar{\beta}}{2 \bar{\beta} + d}} + c (\frac{1}{n})^{1 - \frac{1}{\bar{\beta} + 1}} + c (\frac{1}{n})^{1 - \frac{\bar{\beta}}{2 \bar{\beta} + d}(\frac{1}{\beta_1} + \frac{1}{\beta_2})} \\
& \le c (\frac{1}{n})^{\frac{2 \bar{\beta}}{2 \bar{\beta} + d}} + c (\frac{1}{n})^{\frac{\bar{\beta}}{ \bar{\beta} + 1}} + (\frac{1}{n})^{\frac{2 \bar{\beta} + d -2}{2 \bar{\beta} + d}}
\end{align*}
which is what we wanted, as $d = 2$.
\end{proof}

\section{Proof main results, asynchronous framework}{\label{S: proof asynch}}
In Section \ref{s: asynch} we assume to observe the components of the process $X$ in different moments. {\rev We start by proving the results gathered in Proposition \ref{prop: bound var as}. The proof of Proposition \ref{prop: var due comp continue}, which follows a similar route, can be found in the appendix.}

{\rev 
\subsection{Proof of Proposition \ref{prop: bound var as}}
\begin{proof}
In analogy to the proofs in the synchronous case we introduce
$$k(t , s) := Cov(\prod_{l = 1}^d  {\modarn K_{h_l^*}(x_l - X_{\varphi_{n, l}(t)}^l), \prod_{l = 1}^d K_{h_l^*}(x_l - X_{\varphi_{n, l}(s)}^l)}),$$
such that
\begin{align*}
Var(\hat{\pi}^a_{h^*, T_n}(x)) & = \frac{2}{T_n^2} \int_0^{T_n} \int_0^{t}k(t,s) 1_{s < t} ds dt. 
\end{align*}
We write
\begin{align*}
Var(\hat{\pi}^a_{h^*, T_n}(x)) 
& = \frac{2}{T_n^2} \int_0^{T_n} \int_0^{t}k(t,s) 1_{s < t} \,  \big(1_{|t -s| \le h_1^* h_2^* } \\
& + 1_{ h_1^* h_2^* \le |t -s| \le (\prod_{j \ge 3} h_j^*)^{\frac{2}{d-2}} } + 1_{(\prod_{j \ge 3} h_j^*)^{\frac{2}{d-2}} \le |t -s| \le D  } + 1_{D \le |t -s| \le T_n  } \big)  ds dt \\
& = \sum_{j = 1}^4 \tilde{I}_j, 
\end{align*}
with $h^*$ the rate optimal choice of the bandwidth given by \eqref{eq: optimal bandwidth} {\modarn and $D$ will be specified latter}.
Regarding $\Tilde{I}_1$, acting as in Proposition \ref{prop: bound var discrete dge3} in order to get \eqref{eq: bound k 5.5}, 
{\modarn we have 
\begin{equation*}
\abs{k(t,s)}\le \E \left( \prod_{l = 1}^d K_{h_l^*}^2(x_l - X_{\varphi_{n, l}(t)}^l)\right)^{1/2}	\E \left( \prod_{l = 1}^d K_{h_l^*}^2(x_l - X_{\varphi_{n, l}(s)}^l)\right)^{1/2}.
\end{equation*}
We now state a lemma which will be useful in the sequel and whose proof is postponed to the appendix.
\begin{lemma}\label{L: majo esperance K asynchrone}

1)	We have $\E \left( \prod_{l = 1}^d |K_{h_l}(x_l - X_{\varphi_{n, l}(t)}^l)|\right) \le c$, for some constant $c$ independent of $t$ and $(h_l)_{l=1,\dots,d}$.

2)	We have $\E \left( \prod_{l = 1}^d K_{h_l}^2(x_l - X_{\varphi_{n, l}(t)}^l)\right) \le \frac{c}{\prod_{l=1}^d h_l}$, for some constant $c$ independent of $t$ and $(h_l)_{l=1,\dots,d}$.
\end{lemma}
From application of the second point of Lemma \ref{L: majo esperance K asynchrone},} we obtain $|k(t,s)| \le \frac{c}{\prod_{l = 1}^d h_l^*}.$
Therefore, after the change of variable $t \rightarrow t' := t-s$, we have
\begin{align}{\label{eq: bound I1 asincrono}}
 \tilde{I}_1 & \le \frac{{\rev c}}{T_n^2} \int_0^{T_n} \int_0^{ h_1^* h_2^* \sum_{j = 1}^d |\log h_j^*|} \frac{1}{\prod_{l = 1}^d h_l^*} dt' ds = \frac{c}{T_n} \frac{\sum_{j = 1}^d |\log h_j^*|}{\prod_{l = 3}^d h_l^*},
\end{align}
{\rev where we recall that the constant $c$ may change from line to line.} 
We now study $\tilde{I_2}$, which is the most complicated term we have to deal with. Intuitively, we would like to bound the variance as in the interval $[\delta_1, \delta_2)$ in the proof of Proposition 2 of \cite{Companion}. It relies on a bound for the transition density (as in Proposition 5.1 of \cite{Gobet LAMN} or Lemma 1 of \cite{Companion}). In order to use it we need to know the ordering between the quantities $s$, $t$, $\varphi_{n,1} (s)$, ... , $\varphi_{n,d} (s)$, $\varphi_{n,1} (t)$, ... , $\varphi_{n,d} (t)$. We know that $s < t$, $\varphi_{n,l} (s) \le s $ and $\varphi_{n,l} (t) \le t$ for any $l \in \{ 1, ... , d \}$. Then, we can consider a permutation $w_1$, ... , $w_d$ of $\varphi_{n,1} (s)$, ... , $\varphi_{n,d} (s)$ which is such that $w_1 \le ... \le w_d$. In particular, we denote by $w_1 \le ... \le {\modarn w_d}$ a reordering of $\varphi_{n,l}({\modarn s})$ and $\sigma$ an element of the permutation group on $\{1, ... , d \}$ such that $w_i = \varphi_{n, \sigma(i)}({\modarn s})$ for all $i \in \{ 1, ... , d \}$. In the same way we introduce the permutation $\tilde{w}_1$, ... , $\tilde{w}_d$ of $\varphi_{n,1} (t)$, ... , $\varphi_{n,d} (t)$ which is properly ordered, i.e. $\tilde{w}_1 \le ... \le \tilde{w}_d$. In particular, we introduce
$\tilde{\sigma}$ which is an element of the permutation group such that $\tilde{w}_i = \varphi_{n, \tilde{\sigma}(i)}({\modarn t})$ for all $i \in \{ 1, ... , d \}$. 
As $\Delta_n \le \frac{1}{4} h_1^* h_2^* $ and $|t - s| > h_1^* h_2^* $, we have 
 $$w_1 \le ... \le w_d \le s \le \Tilde{w}_1 \le ... \le \Tilde{w}_d \le t.$$
 We also introduce the vectors, which represent the positions in the instants $w_j$ and $\Tilde{w}_j$. At the instant $w_j$ we have the vector $y^j$, for $j \in \{ 1, ... ,d \}$, while the time $\tilde{w}_j$ are associated to the vectors $\tilde{y}^j$, for $j \in \{ 1 , ... , d \}$.
 %{\modarn We denote $y^j_l$ the $l$-th components of the vector $y^j\in \R^d$.}
  We observe we can write
\begin{align*}
|k(t,s)| & \le |\tilde{k}(t, s)| + |\E[\prod_{l = 1}^d {\modarn K_{h_l^*}(x_l - X_{\varphi_{n, l}(t)}^l)}] | | \E[ \prod_{l = 1}^d {\modarn K_{h_l^*}(x_l - X_{\varphi_{n, l}(s)}^l)}] |\\
& \le |\tilde{k}(t, s)| + c, 
\end{align*}
where 
{\modarn 
	\begin{align*}
		\tilde{k}(t, s) &:= \E[ \prod_{l = 1}^d{\modarn K_{h_l^*}(x_l - X_{\varphi_{n, l}(s)}^l)}  \prod_{l = 1}^d {\modarn K_{h_l^*}(x_l - X_{\varphi_{n, l}(t)}^l)}]
		\\
		&= \E[ \prod_{l = 1}^d K_{h_{\sigma(l)}^*}(x_{\sigma(l)} - X_{w_l}^{\sigma(l)})	\prod_{l = 1}^d K_{h_{\tilde{\sigma}(l)}^*} (x_l - X_{\tilde{w}_l}^{\tilde{\sigma}(l)}) ]
	\end{align*}
}
{\modarn and in the equation above we have used the first point of Lemma \ref{L: majo	esperance K asynchrone}.}
%enleve
%In the equation above we have also used that, thanks to the first point of Lemma \ref{L: technique maj integral K gauss} (stated below and proven in the appendix), there exists a constant $c > 0$ such that
%$$|\E[\prod_{m = 1,2} K(\frac{x_m - X_t^m}{h_m}) \prod_{l = 3}^d K(\frac{x_l - X_{\varphi_{n, l}(t)}^l}{h_l})] | | \E[ \prod_{m = 1,2} K(\frac{x_m - X_s^m}{h_m}) \prod_{l = 3}^d K(\frac{x_l - %X_{\varphi_{n, l}(s)}^l}{h_l})] | \le c.$$} 
We write the expectation in the definition of $\tilde{k}(t,s)$ using the law of the random vector
$(X_{w_1},\dots,X_{w_d},X_{\Tilde{w}_1}, \dots,X_{\Tilde{w}_{d}},)$.
For simplicity, we assume that all the instants appearing in this vector are different $w_1<\dots<w_d<
\Tilde{w}_{1}<\dots<\Tilde{w}_d$ and in turn the law of this vector admits a density as product of the transition density of the process $X$. If we are not in the situation where all the instants are distinct, it is possible to slightly move some values of $w_1,\dots,w_d,\Tilde{w}_{1}, \dots,\Tilde{w}_d$ in order to get different instants, and then conclude by a density argument in order to get the upper bound on $|\tilde{k}(s,t)|$. With these considerations, we can write 
%{\modarn $\tilde{k}(t,s)= \E[ \prod_{l = 1}^d K_{h_{\sigma(l)}}(x_{\sigma(l)} - X_{w_l})
%	\prod_{l = 1}^d K_{h_{\tilde{\sigma}(l)}} (x_l - X_{\tilde{w}_l}^{\tilde{\sigma}(l)}) ]$ as	
%}
\begin{align}{\label{eq: k tilde discrete}}
|\tilde{k}(t,s)| & \le \int_{\R^{d^2}} \prod_{l = 1}^d |K_{h_{\sigma(l)}^*} (x_{\sigma(l)} - y^l_{\sigma(l)})| \int_{\R^{d^2}} \prod_{l = 1}^d |K_{{\modarn h_{\tilde{\sigma}(l)}^*}} (x_{\tilde{\sigma}(l)} - \tilde{y}^l_{\tilde{\sigma}(l)})| p_{w_2 - w_1}(y^1, y^2) p_{w_3 - w_2}(y^2, y^3) \\
& \times ... \times p_{w_d - w_{d - 1}}(y^{d - 1}, y^d) p_{\tilde{w}_{ 1} - w_d}(y^d, \tilde{y}^{ 1}) \times ... \times p_{\tilde{w}_{d} - \tilde{w}_{d-1}}(\tilde{y}^{d-1}, \tilde{y}^{d}) \pi(y^1) dy^1 dy^2... dy^d  d\tilde{y}^{ 1} ...  d\tilde{y}^{d}, \nonumber
\end{align}
{\modarn where we denote $y^j_m$ the $m$-th component of the vector $y^j\in \R^d$.}
 As it is important to remove the contribution of the two smallest {\rev bandwidths}, we need to reorder the components to $\Tilde{y}$. To do that we introduce $\tilde{\sigma}^{-1}$, the inverse of the permutation $\tilde{\sigma}$, and we write
\begin{align*}
|\tilde{k}(t,s)| & \le \int_{\R^{d^2}} \prod_{l = 1}^d |K_{h_{\sigma(l)}^*} (x_{\sigma(l)} - y^l_{\sigma(l)})| \int_{\R^{d^2}} \prod_{l = 1}^d |K_{{\modarn h_{l}^*}} (x_{l} - \tilde{y}^{\tilde{\sigma}^{-1}(l)}_{l})| p_{w_2 - w_1}(y^1, y^2) p_{w_3 - w_2}(y^2, y^3) \\
& \times ... \times p_{w_d - w_{d - 1}}(y^{d - 1}, y^d) p_{\tilde{w}_{ 1} - w_d}(y^d, \tilde{y}^{ 1}) \times ... \times p_{\tilde{w}_{d} - \tilde{w}_{d-1}}(\tilde{y}^{d-1}, \tilde{y}^{d}) \pi(y^1) dy^1 dy^2... dy^d  d\tilde{y}^{ 1} ...  d\tilde{y}^{d}.
\end{align*}
We use the upper bound $\prod_{l=3}^d  \abs{K_{h_{l}^*} (x_{l} - \tilde{y}^{\tilde{\sigma}^{-1}(l)}_{l})} \le C / \prod_{l=3}^d h_l^*$ and we integrate with respect to the variables $\tilde{y}_l^{\sigma^{-1}(l)}$, $l=3,\dots,d$ to get
\begin{multline}\label{E: control kappa var asynch}
|\tilde{k}(t,s)|  \le \frac{C}{\prod_{l=3}^d h_l^*}
\int_{\R^{d^2}}  \prod_{l = 1}^d \abs{K_{h_{\sigma(l)}^*} (x_{\sigma(l)} - y^l_{\sigma(l)})}
\int_{\R^{2d}} \abs{K_{h_1^*}(x_1-\tilde{y}_1^{\tilde{\sigma}^{-1}(1)})
K_{h_2^*}(x_2-\tilde{y}_2^{\tilde{\sigma}^{-1}(2)})} 
\\ \times
p_{w_2 - w_1}(y^1, y^2) p_{w_3 - w_2}(y^2, y^3)  \times ... \times p_{w_d - w_{d - 1}}(y^{d - 1}, y^d)
\\ \times  p_{\tilde{w}_{i_\star} - w_d}(y^d, \tilde{y}^{i_\star}) 
p_{\tilde{w}_{i^\star} - \tilde{w}_{i_\star}}(\tilde{y}^{i_\star}, \tilde{y}^{i^\star}) \pi(y^1)
 dy^1 dy^2... dy^d  d\tilde{y}^{ \tilde{\sigma}^{-1}(1)}d\tilde{y}^{ \tilde{\sigma}^{-1}(2)},
\end{multline}
where $i_\star=\min(\tilde{\sigma}^{-1}(1),\tilde{\sigma}^{-1}(2))$ and 
$i^\star=\max(\tilde{\sigma}^{-1}(1),\tilde{\sigma}^{-1}(2))$.
Using a Gaussian bound on the transition density
 \begin{multline*}% \label{E:def q var asynch}
 	 p_{\tilde{w}_{i_\star} - w_d}(y^d, \tilde{y}^{i_\star}) 
 	p_{\tilde{w}_{i^\star} - \tilde{w}_{i_\star}}(\tilde{y}^{i_\star}, \tilde{y}^{i^\star}) \\
 	\le  \frac{C}{\tilde{w}_{i_\star}-w_d}
 		q\left( (\tilde{y}_j^{\tilde{\sigma}^{-1}(1)})_{j\neq 1},(\tilde{y}_j^{\tilde{\sigma}^{-1}(1)})_{j\neq 2}   \mid y^d,\tilde{y}_1^{\tilde{\sigma}^{-1}(1)}, \tilde{y}_2^{\tilde{\sigma}^{-1}(2)}\right)
 \end{multline*}
where
\begin{multline*}
		q\left( (\tilde{y}_j^{\tilde{\sigma}^{-1}(1)})_{j\neq 1},(\tilde{y}_j^{\tilde{\sigma}^{-1}(1)})_{j\neq 2}   \mid y^d,\tilde{y}_1^{\tilde{\sigma}^{-1}(1)}, \tilde{y}_2^{\tilde{\sigma}^{-1}(2)}\right)
	= \\
	\sqrt{{\modarn \tilde{w}_{i_\star}}-w_d}
	 \prod_{\stackrel{j=1}{j\neq \tilde{\sigma}(i_\star)}}^d
	\frac{e^{-c\frac{( {\modarn \tilde{y}_j^{i_\star}}-y_j^{d} )^2}{ \tilde{w}_{i_\star} -w_d}}}{
	\sqrt{ \tilde{w}_{i_\star} -w_d}}	
\prod_{j=1}^d
\frac{e^{-c\frac{( {\modarn \tilde{y}_j^{i^\star}}-{\modarn \tilde{y}_j^{i^\star}} )^2}{ \tilde{w}_{i^\star} -\tilde{w}_{i_\star}}}}{
	\sqrt{ \tilde{w}_{i^\star} -\tilde{w}_{i^\star}}}.	
\end{multline*}
We now prove that
\begin{multline} \label{E: bound on int q}
\sup_{(y^d,\tilde{y}_1^{\tilde{\sigma}^{-1}(1)},\tilde{y}_2^{\tilde{\sigma}^{-1}(2)}) \in \R^{d + 2}}
	\int_{\R^{2(d-1)}}
		q\left( (\tilde{y}_j^{\tilde{\sigma}^{-1}(1)})_{j\neq 1},(\tilde{y}_j^{\tilde{\sigma}^{-1}(1)})_{j\neq 2}   \mid y^d,\tilde{y}_1^{\tilde{\sigma}^{-1}(1)}, \tilde{y}_2^{\tilde{\sigma}^{-1}(2)}\right)
		\\
\prod_{j=2}^d	d (\tilde{y}_j^{\tilde{\sigma}^{-1}(1)})
\prod_{\stackrel{j=1}{j \neq 2}}^d	d (\tilde{y}_j^{\tilde{\sigma}^{-1}(2)}) \le C.
\end{multline}
To prove \eqref{E: bound on int q}, assume, in order the simplify the notations, that $i_\star=\tilde{\sigma}^{-1}(1)$ and $i^\star=\tilde{\sigma}^{-2}(2)$, as the other case can be proved symmetrically. Then, the integral in the left hand side of \eqref{E: bound on int q}
is
\begin{equation*}
\int_{\R^{2(d-1)}}
	\sqrt{{\modarn\tilde{w}}_{\tilde{\sigma}^{-1}(1)}-w_d}
\prod_{j=2}^d
\frac{e^{-c\frac{( {\modarn \tilde{y}}_j^{\tilde{\sigma}^{-1}(1)}-y_j^{d} )^2}{ \tilde{w}_{\tilde{\sigma}^{-1}(1)} -w_d}}}{
	\sqrt{ \tilde{w}_{\tilde{\sigma}^{-1}(1)} -w_d}}	
\prod_{j=1}^d
\frac{e^{-c\frac{( {\modarn \tilde{y}}_j^{\tilde{\sigma}^{-1}(2)}-{\modarn \tilde{y}}_j^{\tilde{\sigma}^{-1}(1)} )^2}{ \tilde{w}_{\tilde{\sigma}^{-1}(2)} -\tilde{w}_{\tilde{\sigma}^{-1}(1)}}}}{
	\sqrt{ \tilde{w}_{\tilde{\sigma}^{-1}(2)} -\tilde{w}_{\tilde{\sigma}^{-1}(1)}}}
%\\
\prod_{j =2}^d	d (\tilde{y}_j^{\tilde{\sigma}^{-1}(1)})
\prod_{\stackrel{j=1}{j \neq 2}}^d	d (\tilde{y}_j^{\tilde{\sigma}^{-1}(2)}).
\end{equation*}
Integrating with respect to the measures $\prod_{j = 3}^d	d (\tilde{y}_j^{\tilde{\sigma}^{-1}(1)})
\prod_{j =3}^d	d (\tilde{y}_j^{\tilde{\sigma}^{-1}(2)})
$, we get that the last integral is upper bounded by
\begin{equation*}
\int_{\R^2} 
	\sqrt{{\modarn \tilde{w}}_{\tilde{\sigma}^{-1}(1)}-w_d}
\frac{e^{-c\frac{( {\modarn \tilde{y}}_2^{\tilde{\sigma}^{-1}(1)}-y_2^{d} )^2}{ \tilde{w}_{\tilde{\sigma}^{-1}(1)} -w_d}}}{
	\sqrt{ \tilde{w}_{\tilde{\sigma}^{-1}(1)} -w_d}}	
\prod_{j=1}^2
\frac{e^{-c\frac{({\modarn \tilde{y}}_j^{\tilde{\sigma}^{-1}(2)}-{\modarn \tilde{y}}_j^{\tilde{\sigma}^{-1}(1)} )^2}{ \tilde{w}_{\tilde{\sigma}^{-1}(2)} -\tilde{w}_{\tilde{\sigma}^{-1}(1)}}}}{
	\sqrt{ \tilde{w}_{\tilde{\sigma}^{-1}(2)} -\tilde{w}_{\tilde{\sigma}^{-1}(1)}}}
%\\
d \tilde{y}_2^{\tilde{\sigma}^{-1}(1)}
d \tilde{y}_1^{\tilde{\sigma}^{-1}(2)}.
\end{equation*}
Then, integrating with respect to $d\tilde{y}^{\tilde{\sigma}^{-1}(1)}_2$, the convolution of Gaussian kernels yields to the following upper bound for the LHS of \eqref{E: bound on int q},
\begin{equation*}
	\int_{\R} 
	\sqrt{{\modarn\tilde{w}}_{\tilde{\sigma}^{-1}(1)}-w_d}
	\frac{e^{-c\frac{( y_2^{\tilde{\sigma}^{-1}(2)}-y_2^{d} )^2}{ \tilde{w}_{\tilde{\sigma}^{-1}(2)} -w_d}}}{
		\sqrt{ \tilde{w}_{\tilde{\sigma}^{-1}(2)} -w_d}}	
	\frac{e^{-c\frac{( y_1^{\tilde{\sigma}^{-1}(2)}-y_1^{\tilde{\sigma}^{-1}(1)} )^2}{ \tilde{w}_{\tilde{\sigma}^{-1}(2)} -\tilde{w}_{\tilde{\sigma}^{-1}(1)}}}}{
		\sqrt{ \tilde{w}_{\tilde{\sigma}^{-1}(2)} -\tilde{w}_{\tilde{\sigma}^{-1}(1)}}}
	%\\
	d \tilde{y}_1^{\tilde{\sigma}^{-1}(2)}.
\end{equation*}
Using that $\sqrt{ \tilde{w}_{\tilde{\sigma}^{-1}(1)} -w_d} \le \sqrt{ \tilde{w}_{\tilde{\sigma}^{-1}(2)} -w_d}$ and that the first exponential inside the integral above is smaller than $1$, we deduce that \eqref{E: bound on int q} holds true. Using \eqref{E: control kappa var asynch}--\eqref{E: bound on int q} we deduce
\begin{multline*}
	|\tilde{k}(t,s)| 
	\le 
\frac{C}{\prod_{l=3}^d h_l^*} \frac{1}{\tilde{w}_{i_\star}-w_d}
\int_{\R^{d^2}}  \prod_{l = 1}^d \abs{ K_{h_{\sigma(l)}^*} (x_{\sigma(l)} - y^l_{\sigma(l)})}
\int_{\R^{2}}  \abs{K_{h_1^*}(x_1-\tilde{y}_1^{\tilde{\sigma}^{-1}(1)})
		K_{h_2^*}(x_2-\tilde{y}_2^{\tilde{\sigma}^{-1}(2)})} 
	\\ \times
	p_{w_2 - w_1}(y^1, y^2) p_{w_3 - w_2}(y^2, y^3)  \times ... \times p_{w_d - w_{d - 1}}(y^{d - 1}, y^d)
\pi(y^1) dy^1 dy^2... dy^d  d\tilde{y}^{ \tilde{\sigma}^{-1}(1)}_1d\tilde{y}^{ \tilde{\sigma}^{-1}(2)}_2.
\end{multline*}
We now state a lemma which will be useful in the sequel. Its proof can be found in the appendix.

\begin{lemma}\label{L: technique maj integral K gauss}
	1) Let $r \ge 1$ be an integer and $q_1,\dots,q_r \in \{1,\dots,d\}$. Then, we have
	 \begin{equation}  \label{E: integrale en p et K bornee}
	 	\int_{\R^{dr}} \prod_{i=1}^d |K_{h_i}(x_i-u_i^{q_i})|
	 	\prod_{j=1}^{r-1} p_{w_{j+1}-w_{j}}(u^{j},u^{j+1})
	 	du^1\dots du^r \le C,
	 \end{equation}
 for some constant $C$ independent of $(h_i)_i$ and $(w_j)_j$. 
 
 2) The same upper bound holds true if we replace the transition densities $p_{w_{j+1}-w_{j}}$ by Gaussian kernels $g_{(w_{j+1}-w_{j})\lambda_0 Id}$.
\end{lemma}

Using that $\pi$ is bounded, and the {\modarn first} point of Lemma \ref{L: technique maj integral K gauss} with {\modarn $r=d$,} $q_i=\tilde{\sigma}^{-1}(i)$ for $i \in \{1,\dots,d\}$ 
we deduce
{\modarn 
\begin{align*}
|\tilde{k}(t,s)| &\le 
 \frac{c}{\Tilde{w}_{i_\star} - w_d} \frac{1}{\prod_{l \ge 3} h_l^*}
 \int_{\R^{2}}  \abs{K_{h_1^*}(x_1-\tilde{y}_1^{\tilde{\sigma}^{-1}(1)})
 	K_{h_2^*}(x_2-\tilde{y}_2^{\tilde{\sigma}^{-1}(2)})} d\tilde{y}_1^{\tilde{\sigma}^{-1}(1)} d \tilde{y}_2^{\tilde{\sigma}^{-1}(2)}
 		\\
 & \le 
 \frac{c}{\Tilde{w}_{i_\star} - w_d} \frac{1}{\prod_{l \ge 3} h_l^*} \le
\frac{c}{\Tilde{w}_1 - w_d} \frac{1}{\prod_{l \ge 3} h_l^*}
%(1 + \frac{\sqrt{\Delta_n}}{h_2^*}) , 
\end{align*}
}
which implies 
%$$|{k}(t,s)| \le \frac{c}{\Tilde{w}_1 - w_d} \frac{1}{\prod_{l \ge 3} h_l^*}(1 + \frac{\sqrt{\Delta_n}}{h_2^*}) + c .$$
$$|{k}(t,s)| \le \frac{c}{\Tilde{w}_1 - w_d} \frac{1}{\prod_{l \ge 3} h_l^*} + c .$$
In order to bound $\tilde{I_2}$ we also observe that 
\begin{align*}
|t- s| & \le |t - \Tilde{w}_1| + |\Tilde{w}_1 - w_d| + |w_d - s| \\
& \le |\Tilde{w}_1 - w_d|  + 2 \Delta_n \\
& \le |\Tilde{w}_1 - w_d|  + \frac{1}{2} h_1^* h_2^* .
\end{align*}
It implies 
$$\frac{1}{|\Tilde{w}_1 - w_d|} \le \frac{1}{ |t- s| - \frac{1}{2} h_1^* h_2^*}.$$
As a consequence
\begin{align*}
|\tilde{I}_2| & \le \frac{c}{T_n^2} \int_0^{T_n} \int_{0}^{t} \frac{{\modarn 1_{s<t}}}{ |t- s| - \frac{1}{2} h_1^* h_2^*} \frac{1}{\prod_{l \ge 3} h_l^*} 
1_{{\modarn h_1^* h_2^*} \le |t - s | \le {\modarn (\prod_{l \ge 3} h_l^*)^{\frac{2}{d - 2}}}}ds dt.
\end{align*}
By the change of variable $ t-s =: t' $, we obtain
\begin{align*}
|\tilde{I}_2| & \le \frac{c}{T_n} \int_{{\modarn  h_1^* h_2^*} }^{{\modarn  (\prod_{l \ge 3} h_l^*)^{\frac{2}{d - 2}}}} \frac{c}{t' - \frac{1}{2} h_1^* h_2^*} \frac{1}{\prod_{l \ge 3} h_l^*}
% (1 + \frac{\sqrt{\Delta_n}}{h_2^*}) dt' \\
 dt' \\
& \le \frac{c}{T_n} \frac{\sum_{j=1}^2 |\log h_j^*|}{\prod_{l \ge 3} h_l^*}.
%(1 + \frac{\sqrt{\Delta_n}}{h_2^*})
\end{align*}
Regarding $\tilde{I}_3$, it is 
$$\tilde{I}_3 := \frac{1}{T_n^2} \int_0^{T_n} \int_0^{t} k(t,s) 1_{(\prod_{j \ge 3} h_j^*)^{\frac{2}{d-2}} \le |t -s| \le D  } ds dt.$$
We can write $\tilde{k}(t,s)$ as in \eqref{eq: k tilde discrete}.
We use the rough estimation 
$$p_{\tilde{w}_{1} - w_d}(y^d, \tilde{y}^{1}) \le \frac{c}{(\tilde{w}_{1} - w_d)^{\frac{d}{2}}}.$$
We replace it in \eqref{eq: k tilde discrete}, it follows
\begin{align*}
|\tilde{k}(t,s)| & \le \frac{c}{(\tilde{w}_{1} - w_d)^{\frac{d}{2}}} \int_{\R^{d^2}} \prod_{l = 1}^d |K_{h_{\sigma(l)}^*} (x_{\sigma(l)} - y^l_{\sigma(l)})| \int_{\R^{d^2}} \prod_{l = 1}^d |K_{{\modarn h_{\tilde{\sigma}(l)}^*}} (x_{\tilde{\sigma}(l)} - \tilde{y}^l_{\tilde{\sigma}(l)})| \prod_{l = 1} ^{d-1} p_{w_{l + 1} - w_l}(y^l, y^{l + 1}) \\
& \times \prod_{l = 1} ^{d-1} p_{\tilde{w}_{l + 1} - \tilde{w}_{l}}(\tilde{y}^{l}, \tilde{y}^{l + 1}) \pi(y^1) dy^1 dy^2... dy^d  d\tilde{y}^{ 1} ...  d\tilde{y}^{d}. 
\end{align*}
We now apply twice the first point of Lemma \ref{L: technique maj integral K gauss}, having on each integral $r = d$. It provides
\begin{equation}
|k(t,s)| \le \frac{c}{(\tilde{w}_{1} - w_d)^{\frac{d}{2}}} + c.
\label{eq: bound k I3 39.5}
\end{equation}
We now observe it is 
$$|t - s | \le |t - \Tilde{w}_1| + |\tilde{w}_{1} - w_d| + |w_d - s| \le |\tilde{w}_{1} - w_d| + 2 \Delta_n.$$
Hence,
$$ |\tilde{w}_{1} - w_d| \ge |t -s| - 2 \Delta_n \ge |t -s| - \frac{1}{2} h_1^* h_2^* \ge \frac{1}{2} |t - s|.$$
From the change of coordinates $ t' := t-s$ we obtain 
\begin{align*}
\tilde{I}_3 &\le \frac{c}{T_n^2} \int_0^{T_n} \int_{{\modarn  (\prod_{l \ge 3} h_l^*)^{\frac{2}{d-2}}}}^D (\frac{c}{t'^{\frac{d}{2}}} + c) dt' ds \nonumber \\
& \le \frac{c}{T_n} ( (\prod_{l \ge 3} h_l^*)^{\frac{2}{d-2}(1 - \frac{d}{2})} + D) \\
& = \frac{c}{T_n} ((\prod_{l \ge 3} h_l^*)^{-1} + D), \nonumber
\end{align*}
which is the order we wanted. \\
We are left to study $\tilde{I}_4$, which the case where $D \le |t - s| \le T_n$. Here we want to use the fact that the process is exponential $\beta$-mixing. To do that, we use the definition of the covariance. We introduce the notation $K_{h_j^*}(t) := K_{h_j^*}(x - X_t)$. Then we need to study, up to reorder the components, 
$$Cov(K_{h_1^*}(w_1) ...  K_{h_d^*}(w_d), K_{h_1^*}(\tilde{w}_1) ...  K_{h_d^*}(\tilde{w}_d)  ),$$
where $w_1 \le ... \le w_d \le s \le \tilde{w}_1 \le ... \le \tilde{w}_d \le t $, $D \le |t-s| \le T$. We define
$$g(X_{\tilde{w}_1}) := \E[K_{h_1^*}(\tilde{w}_1)  ...  K_{h_d^*}(\tilde{w}_d)|X_{\tilde{w}_1} ].$$
It follows we can write the covariance as 
\begin{align*}
& \E[K_{h_1^*}(w_1) ...  K_{h_d^*}(w_d) K_{h_1^*}(\tilde{w}_1) ... K_{h_d^*}(\tilde{w}_d)] + \\
& - \E[K_{h_1^*}(w_1) ...  K_{h_d^*}(w_d)] \E[K_{h_1^*}(\tilde{w}_1) ... K_{h_d^*}(\tilde{w}_d)] \\
& = \E[K_{h_1^*}(w_1) ...  K_{h_d^*}(w_d) g(X_{\tilde{w}_1}) ] - \E[K_{h_1^*}(\tilde{w}_1) ...  K_{h_d^*}(w_d)] \E[g(X_{\tilde{w}_1})] \\
& = \E[K_{h_1^*}(w_1)  ...  K_{h_d^*}(w_d) (g(X_{\tilde{w}_1}) - \pi(g)) ] \\
& = \E[K_{h_1^*}(w_1) ...  K_{h_d^*}(w_d) (P_{\tilde{w}_1- w_d}g(X_{w_d}) - \pi(g)) ], 
\end{align*}
where we recall that $P_{t} f(x) := \mathbb{E} [f(X_t) | X_0 = x ] = \int_{\mathbb{R}^d} f(y) p_{t} (x,y) dy$ is the transition semigroup of the process $X$.
From {\rew Lemma 7} of \cite{Companion} we get
\begin{align*}%\label{eq: bound k I4 40.5}
|k(t , s)| & \le \prod_{l = 1}^d \left \| K_{h_l^*} \right \|_{\infty} \left \| P_{\tilde{w}_1- w_d}g(X_{w_d}) - \pi(g) \right \|_{L^1} \nonumber \\
& \le \frac{c}{\prod_{l = 1}^d h_l^*} e^{- \rho (\tilde{w}_1- w_d)} \left \|g \right \|_\infty \\
& \le \frac{c}{(\prod_{l = 1}^d h_l^*)^2} e^{- \rho (\tilde{w}_1- w_d)}. \nonumber
\end{align*}
As before, it is clearly $|\tilde{w}_1- w_d| \ge |t - s| - 2 \Delta_n.$ Hence,
\begin{align*}
\tilde{I}_4 & \le \frac{c}{(\prod_{l = 1}^d h_l^*)^2} \frac{1}{T_n^2} \int_0^{T_n} \int_{ D}^{T_n} e^{- \rho s'} e^{\rho \Delta_n} dt ds'\\
& \le \frac{c}{T_n(\prod_{l = 1}^d h_l^*)^2} e^{- \rho D} \nonumber
\end{align*}
Putting all the pieces together, it yields
$$Var(\hat{\pi}^a_{h^*, T_n}(x)) \le \frac{c}{T_n} \frac{\sum_{j = 1}^d |\log h_j^*|}{\prod_{l = 3}^d h_l^*} + \frac{c}{T_n} \frac{\sum_{j = 1}^d |\log h_j^*|}{\prod_{l = 3}^d h_l^*}  + \frac{c}{T_n} \frac{1}{\prod_{l = 3}^d h_l^*} + \frac{D}{T_n} + \frac{c}{T_n (\prod_{l = 1}^d h_l^*)^2} e^{- \rho D}. $$
By choosing $D:= [\max (- \frac{2}{\rho} \log (\prod_{j = 1}^d h_j^*), 1) \land T_n]$ we obtain the wanted result.
\end{proof}

\subsection{Proof of Proposition \ref{prop: bound var as intermediate}}
\begin{proof}
Following the route given by the previous proof we split the variance.
\begin{align*}
Var(\hat{\pi}^a_{\tilde{h}^*, T_n}(x)) 
& = \frac{2}{T_n^2} \int_0^{T_n} \int_0^{t}k(t,s) 1_{s < t} \,  \big(1_{|t -s| \le 3 \Delta_n } \\
& + 1_{ 3 \Delta_n \le |t -s| \le D } + 1_{D \le |t -s| \le T_n  } \big)  ds dt \\
& = \sum_{j = 1}^3 \hat{I}_j.
\end{align*}
We act as on $\Tilde{I}_1$ in the proof of Proposition \ref{prop: bound var as}, obtaining $|k(t,s)| \le \frac{c}{\prod_{l = 1}^d \Tilde{h}^*_l}$ and so $\hat{I}_1 \le \frac{c}{T_n} \frac{\Delta_n}{\prod_{l = 1}^d \Tilde{h}^*_l}$. We bound then $\hat{I}_2$ as $\tilde{I}_3$ in the proof of Proposition \ref{prop: bound var as}. We observe that, in this case, we have $|\Tilde{w}_1 - w_d| \ge |t-s| - 2 \Delta_n$. As $|t-s| \ge 3 \Delta_n$, this implies $|\Tilde{w}_1 - w_d| \ge \frac{1}{3}|t-s|$. The change of coordinates $t':= t-s$, together with \eqref{eq: bound k I3 39.5} yields
$$\hat{I}_2 \le \frac{c}{T_n^2} \int_0^{T_n} \int_{3 \Delta_n}^D (\frac{c}{t'^{\frac{d}{2}}} + c) dt' ds \le \frac{c}{T_n} (\Delta_n^{1 - \frac{d}{2}} + D).$$
We act then on $\hat{I}_3$ as for $\Tilde{I}_4$ above. It follows $$\hat{I}_3 \le \frac{c}{T_n} \frac{1}{(\prod_{l = 1}^d \Tilde{h}^*_l)^2} e^{- \rho D}.$$
We choose $D:= [\max (- \frac{2}{\rho} \log (\prod_{j = 1}^d \Tilde{h}^*_j \, \Delta_n), 1) \land T_n]$. The proof is concluded once one observe that, as $\Delta_n \ge (\prod_{l = 1}^d \Tilde{h}^*_l)^{\frac{2}{d}}$, it is $\Delta_n^{1 - \frac{d}{2}} \le \frac{\Delta_n}{\prod_{l = 1}^d \Tilde{h}^*_l}$.

\end{proof}

}

\subsection{Proof of Proposition \ref{prop: bias asynch}}

\begin{proof}
	From the expression of $\hat{\pi}^a_{h,T_n}(x)$ given in Section \ref{S: Discrete observations} we have
	\begin{equation}\label{E: biais integral temps}
	\E[ \hat{\pi}_{h,{T_n}}^a(x)]=\frac{1}{T_n} \int_0^{T_n} \E \left[ \prod_{l=1}^d K_{h_l}(x_l-X^l_{\varphi_{n,l}(t)}) \right] dt.
\end{equation}
Hence, we focus on $\E \left[ \prod_{l=1}^d K_{h_l}(x_l-X^l_{\varphi_{n,l}(t)}) \right] $for $t \in [0,T_n]$. We denote
by $w_1\le w_2\le \dots \le {\modar w_d}$ a reordering of ${\modar (\varphi_{n,l}(t))_{l=1,\dots,d}}$, and let $\sigma$ an element of the permutation group such that $w_i=\varphi_{n,\sigma(i)}(t)$ for all $i\in \{1,\dots,d\}$. With this notations,
$ \E \left[ \prod_{l=1}^d K_{h_l}(x_l-X^l_{\varphi_{n,l}(t)}) \right]=
\E \left[\prod_{i=1}^d K_{h_{\sigma(i)}}(x_{\sigma(i)}-X^{\sigma(i)}_{w_i}) \right]$ and we now show the following control
\begin{equation}\label{E:maj_principal_biais_sqrt}
\abs{\E \left[\prod_{i=1}^d K_{h_{\sigma(i)}}(x_{\sigma(i)}-X^{\sigma(i)}_{w_i}) \right] - \int_{\mathbb{R}^d} \pi(y)
\prod_{i=1}^d K_{h_i}(x_i-z_i)dz_1 \dots dz_d} \le c \sqrt{{\modarn \Delta_n'}},
\end{equation}
for some constant $c$ independent of $(w_i)_i$ and $(h_i)_i$ {\modarn and where $\Delta_n'$ is defined in \eqref{E: def Delta prime}}.

In the proof of \eqref{E:maj_principal_biais_sqrt} we can assume by a density argument that $w_i<w_{i+1}$ for $i=1,\dots,d-1$. 
In this case, we write
\begin{align} \nonumber
&	\E \left[\prod_{i=1}^d K_{h_{\sigma(i)}}(x_{\sigma(i)}-X^{\sigma(i)}_{w_i}) \right]
	= \int_{\R^{d^2}} \pi(y^1) \prod_{i=1}^d K_{h_{\sigma(i)}}(x_{\sigma(i)}-y_{\sigma(i)}^i) \prod_{i=1}^{d-1} p_{w_{i+1}-w_i}(y^i,y^{i+1}) dy^1 \dots dy^d
\\	\label{E: bias avec intro func xi}
&   \quad\quad\quad\quad
=	\int_{\R^{d}}  \left(\prod_{i=1}^d K_{h_{\sigma(i)}}(x_{\sigma(i)}-y_{\sigma(i)}^i) \right) \xi_{(w_i)_i,\pi} (y_{\sigma(1)}^1,\dots, y_{\sigma(d)}^d) dy_{\sigma(1)}^1 
dy_{\sigma(2)}^2 \dots y_{\sigma(d)}^d,
\end{align}
where for any function $\phi$ we have set
\begin{equation}
	\label{E: def xi}
\xi_{(w_i)_i,\phi} (y_{\sigma(1)}^1,\dots, y_{\sigma(d)}^d)=
\int_{\R^{d(d-1)}} \phi(y^1) \prod_{i=1}^{d-1} p_{w_{i+1}-w_i}(y^i,y^{i+1}) d\widehat{y}^1 \dots d\widehat{y}^d,
\end{equation}
with $\widehat{y}^i=(y^i_j)_{j \in \{1,\dots,d\}\setminus \{ \sigma(i)\}}$. {\rev We now state the following lemma, whose proof can be found in the appendix.}
\begin{lemma}\label{L: biais approx delta xi}
	Suppose that A1 holds and that $a$ and $b$ are $\mathcal{C}^3$ with bounded derivatives. Let $\phi : \mathbb{R}^d \to \mathbb{R}$ be a bounded $\mathcal{C}^1$ function with bounded derivatives and let us denote
$$d_{(w_i)_i, \phi}(y^1_{\sigma(1)},\dots,y^d_{\sigma(d)})=
\xi_{(w_i)_i, \phi}(y^1_{\sigma(1)},\dots,y^d_{\sigma(d)})-
 \phi(y_1^{\sigma^{-1}(1)},\dots,y_d^{\sigma^{-1}(d)}).$$
Then, there exists some {\modar constant} $C$ such that 
%for all $(y^1_{\sigma(1)},\dots,y^d_{\sigma(d)})\in \R^d$, 
we have
\begin{equation*}
	\int_{\R^d}
\left(\prod_{i=1}^d \abs{K_{h_{\sigma(i)}}(x_{\sigma(i)}-y_{\sigma(i)}^i)} \right) 
\abs{ d_{(w_i)_i, \pi}(y^1_{\sigma(1)},\dots,y^d_{\sigma(d)})}
dy_{\sigma(1)}^1\dots dy_{\sigma(d)}^d
 \le C \sqrt{{\modarn \Delta_n'}}.
\end{equation*}
The constant $C$ is independent of $(h_i)_i$.
\end{lemma}

Applying Lemma \ref{L: biais approx delta xi} with $\phi=\pi$, 
%we have $\xi_{(w_i)_i,\pi} (y_{\sigma(1)}^1,\dots, y_{\sigma(d)}^d)=\pi(y_1^{\sigma^{-1}(1)},\dots,y_d^{\sigma^{-1}(d)}) + O_{L^\infty}({\sqrt{\Delta_n}})$. Using this approximation in 
and \eqref{E: bias avec intro func xi}, we deduce
\begin{multline*}
	\E \left[\prod_{i=1}^d K_{h_{\sigma(i)}}(x_{\sigma(i)}-X^{\sigma(i)}_{w_i}) \right]
	\\=
	\int_{\R^{d}}  \left(\prod_{i=1}^d K_{h_{\sigma(i)}}(x_{\sigma(i)}-y_{\sigma(i)}^i) \right) 
	\pi(y_1^{\sigma^{-1}(1)},\dots,y_d^{\sigma^{-1}(d)})
	%\pi(y_{\sigma(1)}^1,\dots, y_{\sigma(d)}^d) 
	dy_{\sigma(1)}^1 \dots dy_{\sigma(d)}^d 
 + O(\sqrt{{\modarn \Delta_n'}}).
\end{multline*}
Changing the notation $z_i=y_i^{\sigma^{-1}(i)}$, which is such that $y_{\sigma(i)}^i=z_{\sigma(i)}$, we get
\begin{align*}
	\E \left[\prod_{i=1}^d K_{h_{\sigma(i)}}(x_{\sigma(i)}-X^{\sigma(i)}_{w_i}) \right]
	&=
\int_{\R^{d}}  \left(\prod_{i=1}^d K_{h_{\sigma(i)}}(x_{\sigma(i)}-z_{\sigma(i)}) \right) 
\pi(z_1,\dots,z_d)
%\pi(y_{\sigma(1)}^1,\dots, y_{\sigma(d)}^d) 
dz_{\sigma(1)} \dots dz_{\sigma(d)} 
+ O(\sqrt{{\modarn \Delta_n'}})
\\&=
\int_{\R^{d}}  \left(\prod_{i=1}^d K_{h_{i}}(x_{i}-z_{i}) \right) 
\pi(z_1,\dots,z_d)
%\pi(y_{\sigma(1)}^1,\dots, y_{\sigma(d)}^d) 
dz_{1} \dots z_{d} 
+ O(\sqrt{{\modarn \Delta_n'}}).
\end{align*}
This implies \eqref{E:maj_principal_biais_sqrt}. Then, recalling \eqref{E: biais integral temps},  and
$ \E \left[ \prod_{l=1}^d K_{h_l}(x_l-X^l_{\varphi_{n,l}(t)}) \right]=\E \left[\prod_{i=1}^d K_{h_{\sigma(i)}}(x_{\sigma(i)}-X^{\sigma(i)}_{w_i}) \right]$
we deduce that
$$	\abs{\E[ \hat{\pi}_{h,{T_n}}^a(x)] - \int_{\R^{d}}  \left(\prod_{i=1}^d K_{h_{i}}(x_{i}-z_{i}) \right) 
	\pi(z_1,\dots,z_d)
	%\pi(y_{\sigma(1)}^1,\dots, y_{\sigma(d)}^d) 
	dz_{1} \dots z_{d}} \le c \sqrt{{\modarn \Delta_n'}}.$$
The proposition follows from the following upper bound on the bias of the synchronous case (see \cite{Chapitre 4})
$$
\abs{
\int_{\R^{d}}  \left(\prod_{i=1}^d K_{h_{i}}(x_{i}-z_{i}) \right) 
\pi(z_1,\dots,z_d) dz_{1} \dots z_{d} - \pi(x_1,\dots,x_d)
} \le c \sum_{i=1}^d h_i^{\beta_i}.
$$
\end{proof}

\subsection{Proof of Theorem \ref{th: conv rate asynch}}
\begin{proof}
The proof is a straightforward consequence of the bias-variance decomposition, Proposition \ref{prop: bound var as} and Proposition \ref{prop: bias asynch}. It is also based on the discussion on the conditions on the discretization step located after the statement of Theorem \ref{th: conv rate asynch}.
\end{proof}

{\rev \subsection{Proof of Theorem \ref{th: conv rate asynch intermediate}}
\begin{proof}
As for the previous theorem, the proof is a direct consequence of the bias-variance decomposition, Proposition \ref{prop: bound var as intermediate} and Proposition \ref{prop: bias asynch}.
\end{proof}}

\begin{appendix}

\section{Appendix}

\subsection{Proof of Proposition \ref{prop: var due comp continue}}
\begin{proof}
In analogy to the previous proofs, we introduce 
$$k(t , s) := Cov(\prod_{m = 1,2} K(\frac{x_m - X_t^m}{h_m^{\rev *}}) \prod_{l = 3}^d K(\frac{x_l - X_{\varphi_{n, l}(t)}^l}{h_l^{\rev *}}), \prod_{m = 1,2} K(\frac{x_m - X_s^m}{h_m^{\rev *}}) \prod_{l = 3}^d K(\frac{x_l - X_{\varphi_{n, l}(s)}^l}{h_l^{\rev *}})),$$
such that
\begin{align*}
Var(\bar{\pi}_{h^*, T_n}(x))
& = \frac{2}{T_n^2} \int_0^{{\rev T_n}} \int_0^{{\rev t}}k(t,s) 1_{s < t} \,  \big(1_{|t -s| \le h_1^* h_2^* \sum_{j = 1}^d |\log h_j^*|} \\
& + 1_{ h_1^* h_2^* \sum_{j = 1}^d |\log h_j^*| \le |t -s| \le (\prod_{j \ge 3} h_j^*)^{\frac{2}{d-2}} } + 1_{(\prod_{j \ge 3} h_j^*)^{\frac{2}{d-2}} \le |t -s| \le D  } + 1_{D \le |t -s| \le T_n  } \big) ds dt \\
& = \sum_{j = 1}^4 I_j, 
\end{align*}
with $h^*$ the rate optimal choice of the bandwidth as in the proof of Theorem 1 of \cite{Companion}, given by \eqref{eq: optimal bandwidth}.
We start considering $I_1$. Acting exactly as in \eqref{eq: bound I1 asincrono} we have
\begin{align*}
 \tilde{I}_1 & \le \frac{c}{T_n} \frac{\sum_{j = 1}^d |\log h_j^*|}{\prod_{l = 3}^d h_l^*}. 
\end{align*}

To study $I_2$ we introduce a similar notation as for the analysis of $\Tilde{I}_2$ for the ordering of $s$, $t$, $\varphi_{n,3} (s)$, ... , $\varphi_{n,d} (s)$, $\varphi_{n,3} (t)$, ... , $\varphi_{n,d} (t)$. In particular, we denote by $w_3 \le ... \le w_d$ a reordering of $\varphi_{n,l}(t)$ and $\sigma$ an element of the permutation group on $\{3, ... , d \}$ such that $w_i = \varphi_{n, \sigma(i)}(t)$ for all $i \in \{ 3, ... , d \}$. In the same way, we 
 introduce
$\tilde{\sigma}$ which is an element of the permutation group such that $\tilde{w}_i = \varphi_{n, \tilde{\sigma}(i)}(s)$ for all $i \in \{ 3, ... , d \}$. %We remark we can see $\Tilde{\sigma} = \sigma \circ \tau$ for some permutation of the same group $\tau$. 
In order to admit the possibility that $\varphi_{n,l}(t) \le s$, and then $\varphi_{n,l}(t) = \varphi_{n,l}(s)$ for some index $l$, we say that $\tilde{w}_j \in \{ w_3, ... , w_d \}$ for $j \le {\rev k}$ and $\Tilde{w}_{{\rev k} + 1} \ge s$.  It follows that 
$$w_3 \le ... \le w_d \le s \le \tilde{w}_{{\rev k} + 1} \le ... \le \tilde{w}_d \le t.$$ Hence we can write, for $l \le {\rev k}$, $\Tilde{w}_l = w_{\tau(l)}$ for some $\tau(l) \in \{3, ..., d\}$.
We now introduce the following vectors, which represent the positions in the instants previously discussed. At the instant $w_j$ we have the vector $y^j$, for $j \in \{ 3, ... ,d \}$. The vector $z$ is instead connected to the time $s$, while $\tilde{z}$ is connected to $t$. At the time $\tilde{w}_j$ we have the vectors $ \tilde{y}^j$, for $j \in \{ {\rev k}+1 , ... , d \}$. We remark that $y^j_l$ is the l-component of the vector $y^j$, which gives the position at the instant $w_j$. We observe that, as $\Tilde{w}_l = w_{\tau(l)}$ for $l \le {\rev k}$ and so we can write, for any $l \le {\rev k}$, $\Tilde{y}_{{\modarn \Tilde{\sigma}} (l)}^l = y_{\sigma(\tau(l))}^{\tau(l)}$. We have
\begin{multline} \label{E: k egal k tilde plus moyenne}
|k(t,s)|  \le |\tilde{k}(t, s)| + 
\\|\E[\prod_{m = 1,2} K(\frac{x_m - X_t^m}{h_m}) \prod_{l = 3}^d K(\frac{x_l - X_{\varphi_{n, l}(t)}^l}{h_l})] | | \E[ \prod_{m = 1,2} K(\frac{x_m - X_s^m}{h_m}) \prod_{l = 3}^d K(\frac{x_l - X_{\varphi_{n, l}(s)}^l}{h_l})] |,
\end{multline}

where 
\begin{equation*}%\label{E: def k tilde asynchrone ctn}
	\tilde{k}(t, s) := \E[\prod_{m = 1,2} K(\frac{x_m - X_t^m}{h_m}) \prod_{l = 3}^d K(\frac{x_l - X_{\varphi_{n, l}(t)}^l}{h_l})  \prod_{m = 1,2} K(\frac{x_m - X_s^m}{h_m}) \prod_{l = 3}^d K(\frac{x_l - X_{\varphi_{n, l}(s)}^l}{h_l})].
	\end{equation*}
 We aim at proving that $|\tilde{k}(t,s)|\le \frac{c}{(t-s) \prod_{l=3}^d
	h_l}$. We have 

%We can write 
\begin{align}{\label{eq: asincrono start 1}}
|\tilde{k}(t,s)| & \le \int_{\R^d} \prod_{m = 1}^2 | K_{h_m^*} (x_m - z_m) | \int_{\R^{d(d - 2)}} |\prod_{l = 3}^d K_{h_{\sigma(l)}^*} (x_{\sigma(l)} - y^l_{\sigma(l)})| \int_{\R^d} \prod_{m = 1}^2 |K_{h_m^*} (x_m - \tilde{z}_m)| \nonumber \\
& \times \prod_{l = 3}^{\rev k} |K_{h_{\sigma(\tau(l))}^*} (x_{\sigma(\tau(l))} - y^{\tau(l)}_{\sigma(\tau(l))})| \int_{\R^{d(d-h)}} \prod_{l = {\rev k} + 1}^d |K_{h_{\tilde{\sigma}(l)}^{\rev *}} (x_{\tilde{\sigma}(l)} - \tilde{y}^l_{\tilde{\sigma}(l)})| \nonumber \\
& \times  p_{w_4 - w_3}(y^3, y^4) p_{w_5 - w_4}(y^4, y^5) \times ... \times p_{w_d - w_{d - 1}}(y^{d - 1}, y^d) p_{s - w_d}(y^d, z) \\
& \times  p_{\tilde{w}_{{\rev k} + 1} - s}(z, \tilde{y}^{{\rev k} + 1}) \times ... \times p_{\tilde{w}_{d} - \tilde{w}_{d-1}}(\tilde{y}^{d-1}, \tilde{y}^{d}) p_{t - \tilde{w}_{d}}( \tilde{y}^{d}, \tilde{z}) \pi(y^3) dz dy^3... dy^d d\tilde{z} d\tilde{y}^{{\rev k} + 1} ...  d\tilde{y}^{d}. \nonumber
\end{align}

We bound
\begin{align*}
& |\prod_{l = 3}^{\rev k} K_{h_{\sigma(\tau(l))}^*} (x_{\sigma(\tau(l))} - y^{\tau(l)}_{\sigma(\tau(l))})\prod_{l = {\rev k} + 1}^d K_{h_{\tilde{\sigma}(l)}^{\rev *}} (x_{\tilde{\sigma}(l)} - \tilde{y}^l_{\tilde{\sigma}(l)})|
%& = |\prod_{l = 3}^{\rev k} K_{h_{\tilde{\sigma}(l)}^*} (x_{\tilde{\sigma}(l)} - \tilde{y}^{l}_{\tilde{\sigma}(l)})\prod_{l = {\rev k} + 1}^d K_{h_{\tilde{\sigma}(l)}^{\rev *}} (x_{\tilde{\sigma}(l)} - \tilde{y}^l_{\tilde{\sigma}(l)})| \\
\le \frac{c}{\prod_{l \ge 3} h_{\tilde{\sigma}(l)}^*} = \frac{c}{\prod_{l \ge 3} h_l^*},
\end{align*}
it follows 
\begin{align}{\label{eq: 24.5}}
|\tilde{k}(t,s)| & \le \frac{1}{\prod_{l \ge 3} h^*_l} \int_{\R^d} {\rev \int_{\R^{d(d - 2)}} } \prod_{m = 1}^2 |K_{h_m^*} (x_m - z_m)| \prod_{l = 3}^d |K_{h_{\sigma(l)}^*} (x_{\sigma(l)} - y^l_{\sigma(l)}) |\int_{\R^d} \prod_{m = 1}^2 |K_{h_m^*} (x_m - \tilde{z}_m)| \nonumber \\
& \times  p_{w_4 - w_3}(y^3, y^4) p_{w_5 - w_4}(y^4, y^5) \times ... \times p_{w_d - w_{d - 1}}(y^{d - 1}, y^d) p_{s - w_d}(y^d, z) \\
& \times p_{t - s}( {y}^{d}, \tilde{z}) \pi(y^3) dz dy^3... dy^d d\tilde{z} . \nonumber
\end{align} 
Using Gaussian upper bounds on the transition density ( Proposition 5.1 in \cite{Gobet LAMN}), we have
%gathered in Lemma 1 of \cite{Companion} we hav
$p_{t - s}( {y}^{d}, \tilde{z}) \le \frac{c}{(t - s)} q_{t - s}(\tilde{z}_3, ... , \tilde{z}_d | \tilde{z}_1, \tilde{z}_2, {y}^d ), $
with
$$q_{t - s}(\tilde{z}_3, ... , \tilde{z}_d | \tilde{z}_1, \tilde{z}_2, {y}^d ) = e^{- \lambda_0 \frac{(\tilde{z}_1 - {y}_1^d)^2}{t - s}} e^{- \lambda_0 \frac{(\tilde{z}_2 - {y}_2^d)^2}{t - s}} \frac{1}{\sqrt{t - s}} e^{- \lambda_0 \frac{(\tilde{z}_3 - {y}_{3}^d)^2}{t - s}} \times ... \times \frac{1}{\sqrt{t - s}} e^{- \lambda_0 \frac{(\tilde{z}_d - {y}_{d}^d)^2}{t - s}}.$$
We observe that 
\begin{equation*}
\sup_{t - s \in (0, 1)} \sup_{{y}^d, \tilde{z}_1, \tilde{z}_2 \in \R^{d + 2}} \int_{\R^{d - 2}} q_{t - s}(\tilde{z}_3, ... , \tilde{z}_d | \tilde{z}_1, \tilde{z}_2, {y}^d ) d\tilde{z}_3 ... d\tilde{z}_d < c,
%\label{eq: q bounded asincrono 2}
\end{equation*}
remarking that $t - s \in (0, 1)$ as $0 \le t-s \le (\prod_{j \ge 3}h_j^*)^{\frac{2}{d-2}} < 1.$
Moreover, we easily bound
$$\int_{\R^2}\prod_{m = 1}^2 K_{h_m^*} (x_m - \tilde{z}_m) d\tilde{z}_1 d\tilde{z}_2 < c.$$
Replacing everything in \eqref{eq: 24.5} we obtain
\begin{align*}%{\label{eq: asincrono k middle 1.5}}
|\tilde{k}(t,s)| & \le \frac{c}{t - s} \frac{1}{\prod_{l \ge 3} h_l^*}  \int_{\R^d} \prod_{m = 1}^2 K_{h_m^*} (x_m - z_m) \prod_{l = 3}^d K_{h_{\sigma(l)}^*} (x_{\sigma(l)} - y^l_{\sigma(l)}) \\
& \times \int_{\R^{d(d - 2)}} \prod_{l = 4}^d p_{w_l - w_{l - 1}}(y^{l - 1}, y^l) p_{s - w_d}(y^d, z) \pi(y^3) dz dy^3... dy^d. \nonumber
\end{align*}
We use the first point of Lemma \ref{L: technique maj integral K gauss} for $r = d-1$. In particular, $u_1, ..., u_{r - 1}, u_r$ are in this case $y^3, ... , y^d, z$ while $q_1 = q_2 = r$ and, for $l \ge 3$, $q_l = \sigma^{-1}(l)$.  
We obtain
$$|\tilde{k}(t,s)| \le \frac{c}{t - s} \frac{1}{\prod_{l \ge 3} h_l^*}. $$
Using again the first point of Lemma \ref{L: technique maj integral K gauss}, there exists a constant $c$ such that
\begin{equation*}
	|\E[\prod_{m = 1,2} K(\frac{x_m - X_u^m}{h_m}) \prod_{l = 3}^d K(\frac{x_l - X_{\varphi_{n, l}(u)}^l}{h_l})] | \le c, \quad \forall u.
\end{equation*}
Recalling \eqref{E: k egal k tilde plus moyenne}, it implies $|{k}(t,s)| \le \frac{c}{t - s} \frac{1}{\prod_{l \ge 3} h_l^*} + c.$
It yields, applying the change of variable $t- s=: t’$,
%\begin{equation*}
%I_2 \le \frac{c}{T_n^2} \int_0^{t}
%\int_0^{T_n} ( \frac{1}{\prod_{l \ge 3} h_l^*} \frac{1}{t - s} + 1) 1_{ h_1^* h_2^* \sum_{j = 1}^d |\log h_j^*| \le |t -s| \le (\prod_{j \ge 3} h_j^*)^{\frac{2}{d-2}} } ds dt.
%\end{equation*}
\begin{align}{\label{eq: I2}}
|I_2| & \le \frac{c}{T_n}   \frac{1}{\prod_{l \ge 3} h_l^*} \int_{h_1^* h_2^* \sum_{j = 1}^d |\log h_j^*|}^{(\prod_{j \ge 3} h_j^*)^{\frac{2}{d-2}}} \frac{1}{t’} dt’ \le \frac{c}{T_n} \frac{\sum_{j=1}^2 |\log h_j^*|}{\prod_{l \ge 3} h_l^*}. 
\end{align}
%$$I_3 := \frac{1}{T_n^2} \int_0^{T_n} \int_0^{T_n} k(t,s) 1_{(\prod_{j \ge 3} h_j^*)^{\frac{2}{d-2}} \le |t -s| \le D  } dt ds.$$
We study now $I_3$. We can write $\tilde{k}(t,s)$ as in \eqref{eq: asincrono start 1}. Now $s$ and $t$ are distant to each other, which implies that $\varphi_{n,l}(s) < \varphi_{n,l}(t)$ for any $l \in \{ 3, ... , d \}$ and so, in particular, the ordering of the quantities previously introduced is the following: 
$$w_3 \le ... \le w_d \le s \le \tilde{w}_3 \le ... \le \tilde{w}_d\le t.$$
This holds true because 
$$|t-s | \ge (\prod_{j \ge 3} h_j^*)^{\frac{2}{d-2}} \qquad \mbox{and } \Delta_n \le (\prod_{j \ge 3} h_j^*)^{\frac{2}{d-2}}.$$
The study of $I_3$ follows the route of the analysis of $\Tilde{I}_3$ in the proof of Proposition \ref{prop: bound var as}. We have 
\begin{align*}
|\tilde{k}(t,s)| & \le \int_{\R^d} \prod_{m = 1}^2 |K_{h_m^*} (x_m - z_m)| \int_{\R^{d(d - 2)}} \prod_{l = 3}^d |K_{h_{\sigma(l)}^*} (x_{\sigma(l)} - y^l_{\sigma(l)})| \int_{\R^d} \prod_{m = 1}^2 |K_{h_m^*} (x_m - \tilde{z}_m)|  \\
& \times \int_{\R^{d(d-2)}} \prod_{l = 3}^d  |K_{h_{\tilde{\sigma}(l)}^*} (x_{\tilde{\sigma}(l)} - \tilde{y}^l_{\tilde{\sigma}(l)})|   p_{w_4 - w_3}(y^3, y^4) p_{w_5 - w_4}(y^4, y^5) \times ... \times p_{w_d - w_{d - 1}}(y^{d - 1}, y^d) p_{s - w_d}(y^d, z) \\
& \times  p_{\tilde{w}_{3} - s}(z, \tilde{y}^{3}) \times ... \times p_{\tilde{w}_{d} - \tilde{w}_{d-1}}(\tilde{y}^{d-1}, \tilde{y}^{d}) p_{t - \tilde{w}_{d}}( \tilde{y}^{d}, \tilde{z}) \pi(y^3) dz dy^3... dy^d d\tilde{z} d\tilde{y}^{3} ...  d\tilde{y}^{d}. \end{align*}
We remark that the largest interval of time above is $\tilde{w}_{3} - s$. We use on it the rough estimation 
$$p_{\tilde{w}_{3} - s}(z, \tilde{y}^{3}) \le \frac{c}{(\tilde{w}_{3} - s)^{\frac{d}{2}}} \prod_{l = 1}^d e^{- \lambda_0 \frac{(z_l - \tilde{y}_{l}^3)^2}{\tilde{w}_3 - s}} \le \frac{c}{(\tilde{w}_{3} - s)^{\frac{d}{2}}}.$$
We apply twice the first point of Lemma \ref{L: technique maj integral K gauss} (having on each integral $r = (d-1)$ as we are considering the integrals in $y^3, ... , y^d, z$ and in $\Tilde{y}^3, ... , \Tilde{y}^d, \Tilde{z}$). It implies 
$$|k(t,s)| \le \frac{c}{(\tilde{w}_{3} - s)^{\frac{d}{2}}} + c.$$
We now observe that it is 
\begin{equation}
|t - s | \le |t - \tilde{w}_3| + |\tilde{w}_3 - s| \le \Delta_n + |\tilde{w}_3 - s| \le \frac{1}{2} (\prod_{l \ge 3} h_l^*)^{\frac{2}{d-2}} + |\tilde{w}_3 - s|.
\label{eq: incremento tempo}
\end{equation}
It follows 
$$|\tilde{w}_3 - s| \ge |t - s | - \frac{1}{2} (\prod_{l \ge 3} h_l^*)^{\frac{2}{d-2}} \ge 
{\modarn \frac{1}{2}|t - s |}.$$
%\frac{1}{2} (\prod_{l \ge 3} h_l^*)^{\frac{2}{d-2}}. $$
Moreover, $|\tilde{w}_3 - s| \le |t - s |  \le D.$
From the change of coordinates  $s \rightarrow s' := t - s$ we obtain 
\begin{align}{\label{eq: bound I3 asincrono}}
I_3 &\le \frac{c}{T_n^2} \int_0^{T_n} \int_{\frac{1}{2} (\prod_{l \ge 3} h_l^*)^{\frac{2}{d-2}}}^D (\frac{c}{s'^{\frac{d}{2}}} + c) ds' dt   = \frac{c}{T_n} ((\prod_{l \ge 3} h_l^*)^{-1} + D),
\end{align}
which is the order we wanted. \\
We are left to study $I_4$, where $D \le |t - s| \le T_n$. Here we act as on $\Tilde{I}_4$ in Proposition \ref{prop: bound var as}. We need to study, up to reorder the components, 
$$Cov(K_{h_1^*}(s) K_{h_2^*}(s) K_{h_3^*}(w_3) ...  K_{h_d^*}(w_d), K_{h_1^*}(t) K_{h_2^*}(t) K_{h_3^*}(\tilde{w}_3) ...  K_{h_d^*}(\tilde{w}_d)  ),$$
where $w_3 \le ... \le w_d \le s \le \tilde{w}_3 \le ... \le \tilde{w}_d \le t $, $D \le |t-s| \le T$. We define
$$g(X_{\tilde{w}_3}) := \E[K_{h_1^*}(t) K_{h_2^*}(t) K_{h_3^*}(\tilde{w}_3) ...  K_{h_d^*}(\tilde{w}_d)|X_{\tilde{w}_3} ].$$
Acting as in the proof of Proposition \ref{prop: bias asynch} we can clearly write the covariance as\\ $\E[K_{h_1^*}(s) K_{h_2^*}(s) K_{h_3^*}(w_3) ...  K_{h_d^*}(w_d) (P_{\tilde{w}_3- s}g(X_s) - \pi(g)) ]$.
%\begin{align*}
%& \E[K_{h_1^*}(s) K_{h_2^*}(s) K_{h_3^*}(w_3) ...  K_{h_d^*}(w_d) K_{h_1^*}(t) K_{h_2^*}(t) K_{h_3^*}(\tilde{w}_3) ... K_{h_d^*}(\tilde{w}_d)] + \\
%& - \E[K_{h_1^*}(s) K_{h_2^*}(s) K_{h_3^*}(w_3) ...  K_{h_d^*}(w_d)] \E[K_{h_1^*}(t) K_{h_2^*}(t) K_{h_3^*}(\tilde{w}_3) ... K_{h_d^*}(\tilde{w}_d)] \\
%& = \E[K_{h_1^*}(s) K_{h_2^*}(s) K_{h_3^*}(w_3) ...  K_{h_d^*}(w_d) (g(X_{\tilde{w}_3}) - \pi(g)) ] \\
%& = \E[K_{h_1^*}(s) K_{h_2^*}(s) K_{h_3^*}(w_3) ...  K_{h_d^*}(w_d) (P_{\tilde{w}_3- s}g(X_s) - \pi(g)) ], 
%\end{align*}
From {\rew Lemma 7} of \cite{Companion} we easily obtain $\left \| P_{\tilde{w}_3- s}g(X_s) - \pi(g) \right \|_{L^1} \le c e^{- \rho (\tilde{w}_3- s)} \left \|g \right \|_\infty.$
Therefore, 
\begin{align*}
|k(t , s)| & \le \prod_{l = 1}^d \left \| K_{h_l^*} \right \|_{\infty} \left \| P_{\tilde{w}_3- s}g(X_s) - \pi(g) \right \|_{L^1} \\
%& \le \frac{c}{\prod_{l = 1}^d h_l^*} e^{- \rho (\tilde{w}_3- s)} \left \|g \right \|_\infty \\
& \le \frac{c}{(\prod_{l = 1}^d h_l^*)^2} e^{- \rho (\tilde{w}_3- s)},
\end{align*}
where we have also used that, from the definition of $g$, it is $\left \|g \right \|_\infty \le \frac{c}{\prod_{l = 1}^d h_l^*}$. Moreover, acting as in \eqref{eq: incremento tempo} and remarking that $(\prod_{l \ge 3} h_l^*)^{\frac{2}{d-2}} \le D$, we easily get
$$|\tilde{w}_3- s| \ge |t - s| - \Delta_n. $$
With the change of variable $s \rightarrow s' := t- s$ we obtain 
\begin{align}{\label{eq: bound I4 asincrono}}
I_4 & \le \frac{c}{(\prod_{l = 1}^d h_l^*)^2} \frac{1}{T_n^2} \int_0^{T_n} \int_{ D}^{T_n} e^{- \rho s'} e^{\rho \Delta_n} dt ds'\le \frac{c}{T_n(\prod_{l = 1}^d h_l^*)^2} e^{- \rho D} 
\end{align}
Putting all the pieces together, using in particular \eqref{eq: bound I1 asincrono}, \eqref{eq: I2}, \eqref{eq: bound I3 asincrono} and \eqref{eq: bound I4 asincrono}, it yields
$$Var(\bar{\pi}_{h^*, T_n}(x)) \le \frac{c}{T_n} \frac{\sum_{j = 1}^d |\log h_j^*|}{\prod_{l = 3}^d h_l^*} + \frac{c}{T_n} \frac{\sum_{j = 1}^d |\log h_j^*|}{\prod_{l = 3}^d h_l^*}  + \frac{c}{T_n} \frac{1}{\prod_{l = 3}^d h_l^*} + \frac{D}{T_n} + \frac{c}{T_n (\prod_{l = 1}^d h_l^*)^2} e^{- \rho D}. $$
By choosing $D:= [\max (- \frac{2}{\rho} \log (\prod_{j = 1}^d h_j), 1) \land T_n]$ we obtain the wanted result.
\end{proof}

%}

\subsection{Proof of Lemma \ref{L: biais approx delta xi}}
\begin{proof}
	Let us denote by $g_\alpha(z)$ the density of a centred Gaussian variable with covariance matrix $\alpha$. We recall the approximation of the diffusion transition density by the Gaussian kernel given in \cite{GobLab08}. Specifying $s=1/N$ and $T=1$ in the notations of the statement of Theorem 3 \cite{GobLab08}, we have for all $s \le 1$,
	\begin{equation}\label{E: control Gobet Labart}
		\abs{p_s(z,z^\prime)- g_{\tilde{a}(z)s}(z-z^\prime)}
		\le C\sqrt{s} g_{{\modarn  s \lambda_0 \text{Id}}}(z-z^\prime) ,
	\end{equation} 	
	where $\lambda_0>0$ and $C>0$ are some constant and $\tilde{a}=a\cdot a^T$. This leads us to introduce a Gaussian approximation of 
	\eqref{E: def xi}
	\begin{equation}
		\label{E: def xi Gaussian}
		\xi_{(w_i)_i,\phi}^{\mathbf{G}} (y_{\sigma(1)}^1,\dots, y_{\sigma(d)}^d)=
		\int_{\R^{d(d-1)}} \phi(y^1) \prod_{i=1}^{d-1} g_{(w_{i+1}-w_i)\tilde{a}(y_i)}(y^{i+1}-y^i) d\widehat{y}^1 \dots d\widehat{y}^d.
	\end{equation}
With this notation, we split $d_{(w_i)_i, \phi}(y^1_{\sigma(1)},\dots,y^d_{\sigma(d)})$ as $
\sum_{l=1}^2
d_{(w_i)_i, \phi}^{(l)}(y^1_{\sigma(1)},\dots,y^d_{\sigma(d)})$,
%+d_{(w_i)_i, \phi}^{(2)}(y^1_{\sigma(1)},\dots,y^d_{\sigma(d)})$,
with
\begin{align} 
	\nonumber
	%\label{E: diff xi xiG}
	d_{(w_i)_i,\phi}^{(1)}(y^1_{\sigma(1)},\dots,y^d_{\sigma(d)})=&
	\xi_{(w_i)_i,\phi} (y_{\sigma(1)}^1,\dots, y_{\sigma(d)}^d)-\xi_{(w_i)_i,\phi}^{\mathbf{G}} (y_{\sigma(1)}^1,\dots, y_{\sigma(d)}^d),
	\\ \label{E: diff xiG phi}
	d_{(w_i)_i,\phi}^{(2)}(y^1_{\sigma(1)},\dots,y^d_{\sigma(d)})=&
\xi_{(w_i)_i,\phi}^{\mathbf{G}} (y_{\sigma(1)}^1,\dots, y_{\sigma(d)}^d)	-\phi(y_1^{\sigma^{-1}(1)},\dots,y_d^{\sigma^{-1}(d)}).	
\end{align}
	The lemma is a consequence of the following upper bound for $l\in\{1,2\}$,
\begin{equation}\label{E: diff xi phi_split}
	\int_{\R^d}
	\left(\prod_{i=1}^d \abs{K_{h_{\sigma(i)}}(x_{\sigma(i)}-y_{\sigma(i)}^i)} \right) 
	\abs{ d_{(w_i)_i, \phi}^{(l)}(y^1_{\sigma(1)},\dots,y^d_{\sigma(d)})}
	dy_{\sigma(1)}^1\dots dy_{\sigma(d)}^d
	\le C \sqrt{\Delta_n}.
\end{equation}
	$\bullet$ We first prove \eqref{E: diff xi phi_split} with $l=1$.
	Comparing  \eqref{E: def xi} {\modar with} \eqref{E: def xi Gaussian}, we can write
	%Using \eqref{E: control Gobet Labart} and \eqref{E: diff xi xiG}, we have
	\begin{multline*}
		d_{(w_i)_i, \phi}^{(1)}(y^1_{\sigma(1)},\dots,y^d_{\sigma(d)})
		= \int_{\R^{d(d-1)}} \phi(y^1)
		\sum_{k=1}^{d-1} \prod_{1\le i<k} p_{w_{i+1}-w_{i}}(y^i,y^{i+1})
		\\
		\times \left[p_{w_{k+1}-w_{k}}(y^k,y^{k+1})- g_{(w_{k+1}-w_{k})
			 \tilde{a}(y^k)}(y^{k+1}-y^k) \right]
		 \prod_{k< i\le d-1} g_{(w_{i+1}-w_{i})\tilde{a}(y^i)}(y^{i+1}-y^i)
		 d\widehat{y}^1 \dots d\widehat{y}^d.
	\end{multline*}
Using \eqref{E: control Gobet Labart} and {\modar a Gaussian upper bound of the transition density,} % under A1 (e.g. see Proposition 5.1 in \cite{Gobet LAMN}), 
we deduce
		\begin{align*}
		\abs{d_{(w_i)_i, \phi}^{(1)}(y^1_{\sigma(1)},\dots,y^d_{\sigma(d)})}
		& \le C \sup_{i=1,\dots,d-1}\sqrt{w_{i+1}-w_i} \times
		\\
		&\quad\quad\quad\quad
		\int_{\R^{d(d-1)}} \abs{\phi(y^1)}
		\prod_{1\le i\le d-1} g_{(w_{i+1}-w_{i})\lambda_0 Id}(y^{i+1}-y^i)
		d\widehat{y}^1 \dots d\widehat{y}^d,
		\\
		& \le C \sqrt{{\modarn \Delta_n'}} \int_{\R^{d(d-1)}} 	\prod_{1\le i\le d-1} g_{(w_{i+1}-w_{i})\lambda_0 Id}(y^{i+1}-y^i)
		d\widehat{y}^1 \dots d\widehat{y}^d,
	\end{align*}
	{\modarn for some constant $\lambda_0>0$,  and where we used that by \eqref{E: def Delta prime}, $w_{i+1}-w_{i}\le \Delta_n'$}. It yields,
	\begin{multline*}
	\int_{\R^d}
\left(\prod_{i=1}^d K_{h_{\sigma(i)}}(x_{\sigma(i)}-y_{\sigma(i)}^i) \right) 
\abs{ d_{(w_i)_i, \phi}^{(1)}(y^1_{\sigma(1)},\dots,y^d_{\sigma(d)})}
dy_{\sigma(1)}^1\dots dy_{\sigma(d)}^d
\\
\le C \sqrt{\Delta_n} \int_{\R^{d^2}}
\prod_{i=1}^d K_{h_{\sigma(i)}}(x_{\sigma(i)}-y_{\sigma(i)}^i) 
	\prod_{1\le i\le d-1} g_{(w_{i+1}-w_{i})\lambda_0 Id}(y^{i+1}-y^i)
d{y}^1 \dots d{y}^d.
	\end{multline*}
From {\rev the second point of Lemma \ref{L: technique maj integral K gauss}} with $r=d$, $u^i=y^i$ for $i=1\dots,d$,
and $q_i=\sigma^{-1}(i)$ we deduce that
 %{\color{red} C EST LA MEME QUE \eqref{eq: asincrono int bounded fine 4} au fait près que l'on integre sur $d$ dates au lieu de $d-2$} 
 the last integral is upper bounded by some constant, and in turn 
 %we deduce
that \eqref{E: diff xi phi_split} holds true for $l=1$.

$\bullet$ We now prove \eqref{E: diff xi phi_split} with $l=2$. In the integral defined by the right hand side of \eqref{E: def xi Gaussian} we make a change of variables,
replacing the variables $(\widehat{y}^1,\dots,\widehat{y}^d)=(y^j_l)_{j \in \{1,\dots,d\}, l \in \{1,\dots,d\}\setminus\{\sigma(j)\}}$ by new integration variables $(z_l^j)_{1\le l \le d, 2 \le j \le d}$ 
defined in the following way. %$z^2,\dots,z^d$  the new integration variables in the following way. 
For $j=d$, we define $z^d=(z^d_1,\dots,z^d_d)$ through the change of variable
$$
y^d_l \rightarrow z^d_l:= \frac{y^d_l-y^{d-1}_l}{\sqrt{w_d-w_{d-1}}}
\text{ for $l \in \{1,\dots,d\}\setminus\{\sigma(d)\}$}, \quad
y^{d-1}_{\sigma(d)} \rightarrow z^d_{\sigma(d)}:= \frac{y^d_{\sigma(d)}-y^{d-1}_{\sigma(d)}}{\sqrt{w_d-w_{d-1}}},
$$
and more generally for $2 \le j \le d$ we define $z^j$ through the formulae
\begin{align*}
y^j_l \rightarrow z^j_l:=& \frac{y^j_l-y^{j-1}_l}{\sqrt{w_j-w_{j-1}}}
\text{ for $l \in \{1,\dots,d\}\setminus\{\sigma(d),\sigma(d-1),\dots,\sigma(j)\}$}, 
\\
y^{j-1}_{l} \rightarrow z^j_{j}:=& \frac{y^j_{l}-y^{j-1}_{l}}{\sqrt{w_j-w_{j-1}}},
	\text{ for $l \in \{\sigma(d),\sigma(d-1),\dots,\sigma(j)\}$}. 
\end{align*}
From these definitions we have
$$
d\widehat{y}^1,\dots,d\widehat{y}^d=\prod_{
	{\modar 	
	%\stackrel{j=1,\dots,d}{l \in \{1,\dots,d\} \setminus \{j\}}
	\substack{j=1,\dots,d \\ l \in \{1,\dots,d\} \setminus \{\sigma(j)\}}
}} dy^j_l=dz^2\dots dz^d 
\prod_{j=2}^d {\modarn (w_j-w_{j-1})^{d/2}}.
$$
{\modar Moreover by construction $z^j=\frac{y^j-y^{j-1}}{\sqrt{w_j-w_{j-1}}}$ for all $j \in \{2,\dots,d\}$.} We deduce that \eqref{E: def xi Gaussian} can be written
{\modar after change of variables}
as
\begin{equation*}
		\xi_{(w_i)_i,\phi}^{\mathbf{G}} (y_{\sigma(1)}^1,\dots, y_{\sigma(d)}^d)=
	\int_{\R^{d(d-1)}} \phi(\widetilde{y}^1) \prod_{j=1}^{d-1} g_{\tilde{a}(\widetilde{y}^j)}(z^{j+1}) dz^2 \dots dz^d,
\end{equation*}
where $\widetilde{y}^j=\widetilde{y}^j(z^2,\dots,z^d;y_{\sigma(1)}^{1},\dots,y_{\sigma(d)}^{d})$ is a notation for the expression of $(y^1,\dots,y^d)$ as a function of the new variables $z^2,\dots,z^d$. They are given by the explicit expression
\begin{equation*}
\widetilde{y}^j_l=\begin{cases}
	y^{\sigma^{-1}(l)}_l+\sum_{u=0}^{j-\sigma^{-1}(l)-1} \sqrt{w_{u+\sigma^{-1}(l)+1}-w_{u+\sigma^{-1}(l)}} z^{u+1+\sigma^{-1}(l)}_l, \quad \text{ 
		if $j > \sigma^{-1}(l)$,}
	\\
	y^{\sigma^{-1}(l)}_l , \quad \text{ if $j = \sigma^{-1}(l)$,}
	\\
	y^{\sigma^{-1}(l)}_l-\sum_{u=0}^{\sigma^{-1}(l)-j-1} \sqrt{w_{u+j+1}-w_{u+j}} z^{u+1+j}_l 
	, \quad\text{ if $j < \sigma^{-1}(l)$.}
\end{cases}	
\end{equation*}
This leads us to introduce the following notation which stresses the dependence upon the time intervals $w_{u+1}-w_{u}$. We set for $s_1,\dots,s_{d-1}>0$,
\begin{equation}\label{E: def y hat s}
\widehat{y}^j_l(s_1,\dots,s_{d-1})=\begin{cases}
	y^{\sigma^{-1}(l)}_l+\sum_{u=0}^{j-\sigma^{-1}(l)-1} s_{u+\sigma^{-1}(l)} z^{u+1+\sigma^{-1}(l)}_l, \quad \text{ 
		if $j > \sigma^{-1}(l)$,}
	\\
	y^{\sigma^{-1}(l)}_l , \quad \text{ if $j = \sigma^{-1}(l)$,}
	\\
	y^{\sigma^{-1}(l)}_l-\sum_{u=0}^{\sigma^{-1}(l)-j-1} s_{u+j} z^{u+1+j}_l 
	, \quad\text{ if $j < \sigma^{-1}(l)$,}
	\end{cases}
\end{equation}
and we define
\begin{equation} \label{E:def xi G s}
	\widehat{\xi^\mathbf{G}} (s_1, \dots, s_{d-1}) =
\int_{\R^{d(d-1)}} \phi(\widehat{y}^1(s_1,\dots,s_{d-1})) \prod_{j=1}^{d-1} g_{\tilde{a}(\widehat{y}^j(s_1,\dots,s_{d-1}))}(z^{j+1}) dz^2 \dots dz^d.
\end{equation}
With these notations, $\widetilde{y}^j_l=\widehat{y}^j_l
(\sqrt{w_2-w_1}, \dots, \sqrt{w_d-w_{d-1}})$ for all $1\le j,l \le d$, 
and
\begin{equation} \label{E: xi dependance inter time} 
	\xi_{(w_i)_i,\phi}^{\mathbf{G}} (y_{\sigma(1)}^1,\dots, y_{\sigma(d)}^d)=	\widehat{\xi^\mathbf{G}} (\sqrt{w_2-w_1}, \dots, \sqrt{w_d-w_{d-1}}).
\end{equation} 
For $(s_1,\dots,s_{d-1})=(0,\dots,0)$ these quantities have simpler expressions. Let us denote $y^\star=(y_1^{\sigma^{-1}(1)}, \dots, y_1^{\sigma^{-1}(d)})$,  and remark that
from \eqref{E: def y hat s}, we have  $\widehat{y}^j(0,\dots,0) = y^\star$, for all $1 \le j \le d$. It follows
\begin{equation}\label{E: xi for s zero}
\widehat{\xi^\mathbf{G}} (0, \dots, 0) 	=
\int_{\R^{d(d-1)}} \phi({y^\star}) \prod_{j=1}^{d-1} g_{\tilde{a}(y^\star)}(z^{j+1}) dz^2 \dots dz^d =
\phi({y^\star})
\end{equation} 
where we have used that $y^\star$ does not depend on the integration variable and that the Gaussian kernel integrates to one.
 We deduce from \eqref{E: xi dependance inter time}--\eqref{E: xi for s zero},
\begin{align} \nonumber
\abs{	\xi_{(w_i)_i,\phi}^{\mathbf{G}} (y_{\sigma(1)}^1,\dots, y_{\sigma(d)}^d)-\phi({y^\star})}
&	=
	\abs{
	\widehat{\xi^\mathbf{G}}(\sqrt{w_2-w_1},\dots,\sqrt{w_{d}-w_{d-1}})-
	\widehat{\xi^\mathbf{G}}(0,\dots,0)}
\\ \label{E: maj xi G - phi}
&\le \sqrt{{\modarn\Delta_n'}} \sum_{j=2}^d \sup_{0\le s_1,\dots,s_{d-1} \le \sqrt{{\modarn\Delta_n'}}} \abs{
\frac{\partial}{\partial s_j} \widehat{\xi^\mathbf{G}}(s_1,\dots,s_{d-1})
}
\end{align}
where we used $\sqrt{w_{j}-w_{j-1}} \le \sqrt{{\modarn\Delta_n'}}$ for all $2\le j \le d$. From the definition \eqref{E: def y hat s}, we have $\abs{\frac{\partial \widehat{y}^i_l}{\partial s_j}} \le c \sum_{u=2}^d \norm{z^u}$. Using that $\phi$ and $\tilde{a}$ are $\mathcal{C}^1$ functions, bounded with bounded derivative, and that $\tilde{a}\ge a^2_{\text{min}} Id$, we deduce from \eqref{E:def xi G s} that
$$
\abs{\frac{\partial \widehat{\xi^\mathbf{G}}(s_2,\dots,s_d)}{\partial s_j}}
\le C \int_{\R^{d(d-1)}} \sum_{u=2}^{d}(1+\norm{z^u}^3) \prod_{j=1}^{d-1} g_{\tilde{a}(\widehat{y}^j(s_1,\dots,s_{d-1}))}(z^{j+1}) dz^2 \dots dz^d.
$$
Used that $\tilde{a}$ is a bounded function under Assumption A1, we deduce that the last integral is bounded independently of $s_1,\dots,s_{d-1}$, and thus 
$\abs{\sup_{s_1,\dots,s_{d-1}}\frac{\partial \widehat{\xi^\mathbf{G}}(s_2,\dots,s_d)}{\partial s_j}} \le C$. In turn, \eqref{E: maj xi G - phi} implies,
\begin{equation*}
	\abs{	\xi_{(w_i)_i,\phi}^{\mathbf{G}} (y_{\sigma(1)}^1,\dots, y_{\sigma(d)}^d)-\phi({y^\star})} \le C \sqrt{{\modarn\Delta_n'}},
\end{equation*}
for a constant $C$ independent of $(w_i)_i$ and $y_{\sigma(1)}^1,\dots, y_{\sigma(d)}^d$. Recalling \eqref{E: diff xiG phi} and the notation
$y^\star=(y_1^{\sigma^{-1}(1)}, \dots, y_1^{\sigma^{-1}(d)})$, this is the upper bound $\abs{d_{(w_i)_i,\phi}^{(2)}
{\modar (y_{\sigma(1)}^1,\dots, y_{\sigma(d)}^d)}}
 \le C \sqrt{{\modarn\Delta_n'}}$, and we deduce \eqref{E: diff xi phi_split} with $l=2$ by integration.
\end{proof}

\subsection{Proof of Lemma \ref{L: technique maj integral K gauss}}

\begin{proof}
	We only prove the first point as the proof of the second  is similar. Using the Gaussian upper bound on the transition density (e.g. see Proposition 5.1 in \cite{Gobet LAMN})
	$$
	p_{w_{j+1}-w_{j}}(u^{j},u^{j+1})\le C \frac{1}{(w_{j+1}-w_j)^{d/2}}e^{-\lambda_0 
		\frac{\abs{u^{j+1}-u^j}^2}{w_{j+1}-w_j}},
	$$
	we deduce that the left hand side of \eqref{E: integrale en p et K bornee} is smaller than 
	\begin{equation*}
		C	\int_{\R^{dr}} \prod_{i=1}^d
			\left[ |K_{h_i}(x_i-u_i^{q_i})| \prod_{j=1}^{r-1} 
			 \frac{1}{(w_{j+1}-w_j)^{1/2}}e^{-\lambda_0 
				\frac{(u^{j+1}_i-u^j_i)^2}{w_{j+1}-w_j}}
			\right] \prod_{i=1}^d \left(\prod_{j=1}^r d u^j_i\right)
	\end{equation*}
which is equal to
\begin{equation}\label{E: prod integrale en dim}
C \prod_{i=1}^d \left( \int_{\R^{r}}  |K_{h_i}(x_i-u_i^{q_i})| \prod_{j=1}^{r-1} 
\frac{1}{(w_{j+1}-w_j)^{1/2}}e^{-\lambda_0 
	\frac{(u^{j+1}_i-u^j_i)^2}{w_{j+1}-w_j}}
du^1_i \dots du^r_i	\right).
\end{equation}
It is sufficient to show that for all $i \in \{1,\dots,d\}$ the integrals in the product \eqref{E: prod integrale en dim} are smaller than some constant independent of $(h_i)_i$ and $(w_j)_j$.
We successively integrate in $u^r_i,u^{r-1}_i\dots,u^{q_i+1}_i$, using the change of variables 
$u^j_i \rightarrow z^j_i:=\frac{u^{j}_i-u^{j-1}_i}{\sqrt{w_{j+1}-w_j}}
$ 
for $j=r,r-1,\dots,q_i+1$. Next, we successively integrate in $u^1_i,\dots,u_i^{q_i-1}$, using the change of variable $u^j_i \rightarrow z^j_i:=\frac{u^{j+1}_i-u^{j}_i}{\sqrt{w_{j+1}-w_j}}$. We deduce that the integrals appearing in \eqref{E: prod integrale en dim} are upper bounded by $C \int_\R  K_{h_i}(x_i-u_i^{q_i}) d u_i^{q_i} \le C$. This proves the lemma.
\end{proof}

{\modarn
\subsection{Proof of Lemma \ref{L: majo esperance K asynchrone}}

\begin{proof}
	We start by the proof of the first point. We reorder $(\varphi_{n,l}(t))_{l=1,\dots,d}$ as $w_1 \le \dots \le w_d$ and denote by $\sigma$ a permutation of $\{1,\dots,d\}$ such that
	$\varphi_{n,\sigma(l)}(t)=w_l$ for all $l \in \{1,\dots,d\}$. By a density argument we can assume that the $w_l$ are all 	distinct. Then, we can write
	\begin{equation*}
		\E\left[\prod_{l=1}^d |K_{h_l}(x_l-X^l_{\varphi_{n,l}(t)})|\right]=\int_{\mathbb{R}^{d^2}} \prod_{l=1}^d
		|K_{h_l}(x_l-y_l^{\sigma^{-1}(l)})| \prod_{l=1}^{d-1} p_{w_{l+1}-w_l}(y^l,y^{l+1})dy^1\dots dy^d.
	\end{equation*}
	Now the first part of the lemma is a consequence of the first point of Lemma \ref{L: technique maj integral K gauss}, with $r=d$ and $q_l=\sigma^{-1}(l)$ for all $l \in \{1,\dots,d\}$.
	
	The second point of the lemma is obtained as a consequence of the first point, after remarking that we can write $K_{h_l}(\cdot)^2 = \frac{1}{h_l} (\frac{1}{h_l}K^2(\frac{\cdot}{h_l}))$, and applying the first point with the function $K^2$ instead of $|K|$.
\end{proof}
}

\end{appendix}

\end{document}